%% file: main.tex
\documentclass{amsart}

\usepackage{amssymb}

\usepackage[cmtip,all]{xy}
\usepackage[utf8]{inputenc}
\usepackage{tikz}
\usepackage{tikz-cd}
\usepackage{todonotes}
\usepackage{rotating}
\usepackage{hyperref}
\usepackage{centernot}
\usepackage{verbatim}
\usepackage[justification=centering]{caption}
\usetikzlibrary{matrix,shapes,arrows,arrows.meta,calc,topaths,intersections,hobby,positioning,decorations.pathreplacing,decorations.pathmorphing,fit,patterns}
\tikzset{
asdstyle/.style={blue,thick},
righthairs/.style={postaction={decorate,draw,decoration={border,amplitude=\hramp,segment length=\hrlen,angle=-90,pre=moveto,pre length=\hrlen/2}}},
lefthairs/.style={postaction={decorate,draw,decoration={border,amplitude=\hramp,segment length=\hrlen,angle=90,pre=moveto,pre length=\hrlen/2}}},
righthairsnogap/.style={postaction={decorate,draw,decoration={border,amplitude=\hramp,segment length=\hrlen,angle=-90}}},
lefthairsnogap/.style={postaction={decorate,draw,decoration={border,amplitude=\hramp,segment length=\hrlen,angle=90}}},
graphstyle/.style={thick},
arrowstyle/.style={thick,decorate,decoration={snake,amplitude=1.7,segment length=10pt,post length=.5mm,pre length=0}},
genmapstyle/.style={thick,-stealth'},
arrhdstyle/.style={thick},
exceptarcstyle/.style={red, ultra thick},
dualquiverstyle/.style={thick,->},
|/.tip={Bar[width=.8ex]}
}

\usepackage{pict2e}
\makeatletter
\newcommand{\adjunction}[4]{%
  #1\colon #2%
  \mathrel{\vcenter{%
    \offinterlineskip\m@th
    \ialign{%
      \hfil$##$\hfil\cr
      \longrightharpoonup\cr
      \noalign{\kern-.3ex}
      \smallbot\cr
      \longleftharpoondown\cr
    }%
  }}%
  #3 \noloc #4%
}
\newcommand{\longrightharpoonup}{\relbar\joinrel\rightharpoonup}
\newcommand{\longleftharpoondown}{\leftharpoondown\joinrel\relbar}
\newcommand\noloc{%
  \nobreak
  \mspace{6mu plus 1mu}
  {:}
  \nonscript\mkern-\thinmuskip
  \mathpunct{}
  \mspace{2mu}
}
\newcommand{\smallbot}{%
  \begingroup\setlength\unitlength{.15em}%
  \begin{picture}(1,1)
  \roundcap
  \polyline(0,0)(1,0)
  \polyline(0.5,0)(0.5,1)
  \end{picture}%
  \endgroup
}
\makeatother

\newcommand{\po}{\ar@{}[dr]|{\text{\pigpenfont R}}}
\newcommand{\pb}{\ar@{}[dr]|{\text{\pigpenfont J}}}

\newcommand*{\hrlen}{10}
\newcommand*{\hramp}{3}

\newcommand{\yo}{\text{\usefont{U}{min}{m}{n}\symbol{'110}}}
\DeclareFontFamily{U}{min}{}
\DeclareFontShape{U}{min}{m}{n}{<-> dmjhira}{}

\theoremstyle{definition}
\newtheorem{theorem}{Theorem}[section]
\newtheorem{lemma}[theorem]{Lemma}
\newtheorem{prop}[theorem]{Proposition}
\newtheorem{cor}[theorem]{Corollary}
\newtheorem{definition}[theorem]{Definition}
\newtheorem{example}[theorem]{Example}

\newtheorem{notation}[theorem]{Assumption}
\newtheorem{construction}[theorem]{Construction}
\theoremstyle{remark}
\newtheorem{remark}[theorem]{Remark}

\numberwithin{equation}{section}

\newcommand{\A}{\mathbb{A}}
\newcommand{\Gm}{\mathbb{G}_m}

\newcommand{\ind}{\operatorname{Ind}}
\newcommand{\interior}{\operatorname{int}}
\newcommand{\coh}{\operatorname{Coh}}
\newcommand{\sh}{\operatorname{Sh}}
\newcommand{\psh}{\operatorname{PSh}}
\newcommand{\ush}{\operatorname{\mu Sh}}
\newcommand{\uhom}{\operatorname{\mu hom}}

\newcommand{\fun}{\operatorname{Fun}}
\newcommand{\map}{\operatorname{Fun}}

\newcommand{\DGCat}{\operatorname{DGCat}_k}
\newcommand{\dgCat}{\operatorname{dgCat}_k}
\newcommand{\oo}{\infty}
\newcommand{\ustk}{\operatorname{\mu stk}}
\newcommand{\ssupp}{\operatorname{SS}}
\newcommand{\open}{\operatorname{Open}}

\newcommand{\op}{\mathrm{op}}
\newcommand{\Z}{\mathbb{Z}}
\newcommand{\R}{\mathbb{R}}
\newcommand{\C}{\mathbb{C}}
\newcommand{\cat}{\operatorname{Cat}}
\newcommand{\sset}{\operatorname{sSet}}
\newcommand{\pic}{\operatorname{Pic}}
\newcommand{\abs}[1]{\lvert #1 \rvert}
\usepackage{rotating}
\newcommand*{\isoarrow}[1]{\arrow[#1,"\rotatebox{90}{\(\sim\)}"]}

\newcommand{\atwo}{\bullet\rightarrow\bullet}
\newcommand{\athree}{\bullet\leftarrow\bullet\rightarrow\bullet}
\newcommand{\spec}{\operatorname{Spec}}
\newcommand{\rhom}{\operatorname{RHom}}
\newcommand{\Ac}{\mathcal{A}}
\newcommand{\Oc}{\mathcal{O}}
\newcommand{\Ec}{\mathcal{E}}
\newcommand{\Fc}{\mathcal{F}}
\newcommand{\Xc}{\mathcal{X}}

\newcommand{\Cc}{\mathcal{C}}
\newcommand{\Dc}{\mathcal{D}}
\newcommand{\Lc}{\mathcal{L}}
\renewcommand{\Mc}{\mathcal{M}}
\newcommand{\Sym}{\mathcal{S}ym}
\newcommand{\colim}{\operatornamewithlimits{colim}}
\newcommand{\Spec}{\operatorname{\mathcal{S}pec}}
\newcommand{\coker}{\operatorname{coker}}
\newcommand{\Cat}{\mathcal{C}at}
\newcommand{\Mod}{\operatorname{\mathbf{Mod}}}
\newcommand{\id}{\operatorname{id}}
\begin{document}

% \title[short text for running head]{full title}
\title{Coherent-Constructible Correspondence for Toric Fibrations}

%    Only \author and \address are required; other information is
%    optional.  Remove any unused author tags.

%    author one information
% \author[short version for running head]{name for top of paper}
\author{Yuxuan Hu}
\address{Department of Mathematics, Northwestern University, 2033 Sheridan Road, Evanston, IL 60208, USA}
\curraddr{}
\email{yhu@math.northwestern.edu}
\thanks{}

%    author two information
\author{Pyongwon Suh}
\address{Department of Mathematics,
 Lady Shaw Building,
The Chinese University of Hong Kong,
Shatin, N.T., Hong Kong}
\curraddr{}
\email{pyongwonsuh@cuhk.edu.hk}
\thanks{}

%    \subjclass is required.
\subjclass[2010]{}

\date{}

\dedicatory{}

%    Abstract is required.
\begin{abstract}
Let $\Sigma$ be a fan inside the lattice $\mathbb{Z}^n$,
and $\Ec:\Z^n \rightarrow \pic{S}$ be a map of abelian groups.
We introduce the notion of a principal
toric fibration $\Xc_{\Sigma, \Ec}$ over the base scheme $S$,
relativizing the usual toric construction for $\Sigma$.
We show that the category of ind-coherent
sheaves on such a fibration is equivalent to the
global section of the Kashiwara-Schapira
stack twisted by a certain local system
of categories with stalk $\ind\coh S$.
It is a simultaneous generalization of the work of Harder-Katzarkov \cite{harder2019perverse} and of Kuwagaki \cite{kuwagaki20},
and should be seen as a family-version of the coherent-constructible correspondence \cite{fang2011categorification}.
\end{abstract}

\maketitle
%    Text of article.
\tableofcontents
\input{sections/introduction}
\input{sections/notations}
% \input{sections/intro}
% \input{sections/local-case}
\input{sections/affine}
\input{sections/affine-stacky}
\input{sections/descent}
\input{sections/quiver}

%    Bibliographies can be prepared with BibTeX using amsplain,
%    amsalpha, or (for "historical" overviews) natbib style.
\bibliographystyle{amsalpha}
%    Insert the bibliography data here.
\bibliography{ref.bib}
\end{document}

%% file: sections/Introduction.tex
\section{Introduction}
\subsection{Mirror symmetry and coherent-constructible correspondence}
Homological Mirror Symmetry (HMS) is an equivalence between two categories of seemingly different geometric origins. One of these categories is the derived category of coherent sheaves on a complex variety $X$. The other is the Fukaya-type category of another seemingly unrelated variety $\check{X}$ or a potential function $W:Y\to \C$, whose objects are roughly Lagrangian submanifolds and whose morphisms are their intersections. The composition of morphisms in the latter category is defined by counting pseudo-holomorphic curves.

On the other hand, it has been known that the Fukaya category of an exact symplectic manifold is sheaf theoretic in nature. For example, the Nadler-Zaslow theorem \cite{nadlerzaslow, nadler2009} states that there is an equivalence between the infinitesimal Fukaya category of the cotangent bundle of a real analytic manifold $Z$ and the category of constructible sheaves on $Z$. Kontsevich \cite{kontsevich2009} conjectured that for a Weinstein manifold, there is a (co)sheaf of categories on its Lagrangian skeleton whose global section is equivalent to the Fukaya category of the Weinstein manifold. Nadler's work proposed that this (co)sheaf of categories should be modeled by the category of microlocal sheaves or a wrapped version of them \cite{nadler2016}. This is now a theorem for many cases due to the work of Ganatra-Pardon-Shende \cite{ganatrapardonshende1, ganatrapardonshende2, ganatrapardonshende3}.

The Coherent-Constructible Correspondence (CCC) is a version of HMS for toric varieties which intertwines with the sheaf theoretic description of the Fukaya category. Originally conjectured by Bondal \cite{bondal}, the equivariant version was proved by Fang-Liu-Treumann-Zaslow \cite{fang2011categorification}. Some special cases of the non-equivariant version were proved in \cite{treumann2010remarks, scherotzkesibilla, kuwagaki2017}, and now it is a theorem which can be proved for toric stacks \cite{kuwagaki20, zhou19}, which we describe as follows.

Starting from a fan datum $\Sigma$ in the lattice $N$,
one can construct a $n$-dimensional toric variety $X_\Sigma$ as well
as a singular Lagrangian skeleton
\[
\Lambda_\Sigma=\bigcup_{\sigma \in \Sigma} (\sigma^{\perp})\times (-\sigma)  \subseteq T^*M_T.
\]
Here $M_T = N^\vee \otimes_\Z \R \simeq T^n$ is the dual torus. 
We can then study the category $\sh_{\Lambda_{\Sigma}}(T^n)$ of constructible sheaves
with singular support contained in $\Lambda_\Sigma$.
\begin{theorem}[CCC \cite{kuwagaki20}]
There is an equivalence of $\infty$-categories
\[\coh(X_\Sigma) \simeq \sh^w_{\Lambda_\Sigma}(T^n),\]
where $\coh(X_\Sigma)$ is the bounded derived category of coherent sheaves on $X_\Sigma$, and $\sh^w_{\Lambda}(T^n)$
is the category of wrapped sheaves\footnote{The category of wrapped sheaves is the full subcategory of compact objects in the category of unbounded sheaves
(with singular supports contained in $\Lambda$).}.
\end{theorem}

We remark how this is related to
homological mirror symmetry in concrete terms.\footnote{
For the sake of simplicity,
we consider everything up to Morita equivalance and don't distinguish
a small category and its ind-completion.
This is possible because the ind-completion functor
gives an equivalence between the category
of small dg-categories localized at Morita equivalences and the category
of cocomplete dg-categories.
}
By work of \cite{nadlerzaslow, nadler2009, ganatrapardonshende3},
the sheaf category $\sh_\Lambda(M)$ is
equivalent to the partially wrapped
Fukaya category $\mathcal{W}(T^*M, \Lambda)$.
For a toric variety, one can also
consider its mirror Landau-Ginzburg model
$W:(\C^\times)^n \rightarrow \C$,
which can be written explicitly using the
fan data.
The Floer-theoretic category we want
to consider in this case is a variant $\mathcal{W}((\C^\times)^n, W^{-1}(\infty))$ of
the Fukaya-Seidel category.
Note that symplectically $(\C^\times)^n \simeq T^*T^n$, and in \cite{rstz14, zhou20}, it is shown that under this identification
the hypersurface $W^{-1}(\infty)$
admits a Liouville Lagrangian skeleton
given by (the contact infinity part of) the FLTZ skeleton $\Lambda_\Sigma$.
Finally, note that, by \cite{ganatrapardonshende1}, replacing the stop
$W^{-1}(\infty)$ with its Liouville skeleton 
doesn't change the partially wrapped Fukaya category.
This way, we see that the category
$\sh_{\Lambda_\Sigma}(T^n)\simeq \mathcal{W}((\C^\times), W^{-1}(\infty))$ is really
a sheaf-theoretic incarnation of the
Fukaya-Seidel category of the mirror
Landau-Ginzburg model to $X_\Sigma$.
\subsection{CCC in the relative case}
In the microlocal sheaf theory of Kashiwara-Schapira \cite{kashiwaraschapira}, there is a sheaf of categories  $\mu sh_\Lambda$ on $T^*M$ called the Kashiwara-Schapira stack whose category of global sections recovers the category of constructible sheaves  \cite[Proposition 3.5]{shendetreumannwilliams}:

\[\mathrm{Sh}_\Lambda(M) \simeq \Gamma(\Lambda, \ush_\Lambda).\]
By taking ind-completion, the coherent-constructible correspondence can be stated as
\[\ind\coh(X_\Sigma) \simeq \Gamma(\Lambda_{\Sigma}, \ush_{\Lambda_{\Sigma}}).\]
If we consider the product of the toric variety with
another variety $S$, by the K\"unneth formula, we have
\[
\ind\coh(X_\Sigma\times S) \simeq \Gamma(\Lambda_\Sigma, \ush_{\Lambda_\Sigma})\otimes_k \ind\coh(S)\simeq \Gamma(\Lambda_\Sigma, \ush_{\Lambda_\Sigma}\otimes_k \underline{\ind\coh(S)}),
\]
where $\underline{\ind\coh(S)}$ is the constant
sheaf of categories with stalk $\ind\coh(S)$.
Here the tensor products are tensor products of (sheaves of) $k$-linear
stable $\oo$-categories, see \ref{subsec:tensor-products}.

The topic of this paper is to generalize this equivalence from the case of a Cartesian product $X_\Sigma \times S$ to that of a \emph{principal toric fibration} (see \ref{subsection:fibration}) over $S$.
Given any fan datum $\Sigma$ and sequence of line bundles $\mathcal{E}=(\Lc_1, \dots, \Lc_n)$ one can define such a principal toric fibration, whose fiber is the toric variety $X_\Sigma$ associated to the fan $\Sigma$.
The line bundles $\mathcal{E}$ give rise to a local system of categories $\Phi(\mathcal{E})$ on $T^n$ whose monodromy along the $i$-th generator of $\pi_1(T^n)\simeq \Z^n$ is given by $-\otimes \Lc_i^\vee$ ( see Remark \ref{ikari-altenative-defn}, and Sections \ref{subsection:fibration}
\ref{subsection:local-system} for detailed discussions.)

In the setting of principal toric fibrations, we will show that the equivalence continues to hold provided the constant sheaf is replaced with the local system of categories $\Phi(\Ec)$ (see definition \ref{def:local-system})
with stalk $\ind\coh(S)$.
More precisely, our main theorem is the following:
% We will prove CCC for toric fiber bundles. By a toric fiber bundle on a variety $Y$, we mean a fiber bundle $\mathcal{X}_{\Sigma}$ on $Y$ whose fiber is a toric variety $X_\Sigma$ and the structure group is $(\C^*)^n$ which lies in and acts naturally on $X_\Sigma$.
\begin{theorem}[Theorem \ref{maintheorem}]
Let $\mathcal{X}_{\Sigma, \Ec}$ be a (stacky) principal toric fibration over $S$ constructed from the fan $\Sigma$ and the sequence of line bundles $\Ec$ (hence whose fiber is the toric variety $X_\Sigma$). Then there is an equivalence.

\[\ind\coh (\mathcal{X}_{\Sigma,\Ec}) \simeq \Gamma(\Lambda_\Sigma, \ush_{\Lambda_\Sigma} \otimes_k \Phi(\Ec)). \]

% where

% \begin{itemize}

%     \item $\ind\coh(\mathcal{X}_\Sigma)$ is the $\infty$-category of ind-coherent sheaves on $\mathcal{X}_\Sigma$.
    
%     \item $\Lambda_\Sigma$ is the FLTZ skeleton associated to $\Sigma$.

%     \item $\ush_{\Lambda_\Sigma}^\diamond$ is the Kashiwara-Schapira stack associated to $\Lambda_{\Sigma}$.
    
%     \item $\Phi(\mathcal{X}_\Sigma)$ is the local system of categories on $T^*T^n$ hence on $\Lambda_\Sigma$ whose generic stalk is $\ind\coh(Y)$ and monodromies are $-\otimes \mathcal{M}_i$: recall that associated to $Z$ is a $(\C^*)^n$-bundle hence a direct sum of line bundles $\mathcal{M}_1 \oplus \varietyoplus \mathcal{M}_n$.
% \end{itemize} 
\end{theorem}

\subsection{Perverse schobers on surfaces}
The idea of CCC for toric fibrations dates back to the work of Harder-Katzarkov \cite{harder2019perverse}.
Let $\mathcal{O}$ be the structure sheaf of a variety $X$ and $\mathcal{L}$ be a line bundle on $X$. 
% Abusing notations, we write $\Oc$ and $\Lc$ for their total space.
Then $\mathbb{P}_X(\Oc \oplus \Lc)$, the projectivization of the rank $2$ vector bundle $\Oc \oplus \Lc$, is a toric $\mathbb{P}^1$-bundle on $X$. \footnote{In fact, any toric $\mathbb{P}^1$-bundle on $X$ is of this form} Their theorem can be restated as follows using our notation.
\begin{theorem}[\cite{harder2019perverse} Theorem 4.9]
There is a equivalence of $k$-linear stable
$\oo$-categories.
\[\ind\coh(\mathbb{P}_X(\Oc \oplus \Lc)) \simeq \Gamma(T^*T^1, \ush_{\Lambda_{\Sigma_{\mathbb{P}^{1}}}}\otimes_k \Phi(\mathbb{P}_X(\Oc \oplus \Lc)))\]
\end{theorem}

In fact, Theorem 4.9 as stated in \cite{harder2019perverse} is different from the above. They used their own definition of a perverse sheaf of categories on a surface which originated from the idea of perverse schobers due to Kapranov-Schechtman \cite{kapranovschechtman}.
In addition to their theorem, it was conjectured that there should be a version of CCC for toric fiber bundles.
Our main theorem is a realization of their conjecture.

Harder-Katzarkov did not make a precise statement of their conjecture.
%One possible reason is that they used their own definition of perverse sheaf of categories on a surface which originated from the idea of perverse schobers due to Kapranov-Schechtman \cite{kapranovschechtman}, which is quiver theoretic in nature, whose higher-dimensional analog defined only in special cases as of now.
However, their viewpoint provides a motivation for a quiver theoretic description of the constructible-side category $\Gamma(T^*T^n, \ush_{\Lambda_\Sigma} \otimes \Phi(\mathcal{X}_\Sigma))$. Such a quiver-theoretic description provides a description of the derived category of a toric fiber bundle in terms of the derived categories of its base and fiber. This might be of independent interest. As an application of this approach, in section \ref{quiver}, we give a quiver theoretic proof of the main theorem for the case of toric $\mathbb{P}^n$-bundles.

\subsection{On Fukaya categories}
We remark on what we should expect for a Floer-theoretic counterpart of our result.
\footnote{We warn the reader that the statements in this section are speculative and only ``morally true"
per se, given that the equivalence between different flavors of Fukaya-Seidel category is yet to be established to our knowledge, among
other subtleties.
However we believe one with sufficient
expertise in symplectic
geometry are able to carry out
complete proofs of these claims
once those techinal issues are resolved.
}
This has also been dicussed in \cite[remark 4.8]{harder2019perverse} in the $\mathbb{P}^1$ case.
We will work with Liouville sectors as in \cite{ganatrapardonshende1} on which
wrapped Fukaya categories have better functoriality.
Note that the partially wrapped Fukaya category
with stop a hypersurface is the same as the wrapped
Fukaya category of the Liouville sector obtained by
removing (a neighborhood of) the stop.

From this point of view, denote by $X_\Sigma^\vee$ the sector obtained from the mirror Landau-Ginzburg model
to a toric variety $X_\Sigma$, and suppose we have
a mirror pair $(S, S^\vee)$ so that \[
\coh(S)\simeq \mathcal{W}(S^\vee).
\]
Further, suppose there is a fibration $\Xc_\Sigma^\vee$ of Liouville
sectors over $X_\Sigma^\vee$ with fiber $S^\vee$ in a suitable sense,
so that the symplectic monodromy acts
via the Picard group $\operatorname{Pic}(S)$ when transferred
to the coherent side via the equivalence.
From this monodromy data,
one can construct a principal toric fibration $\Xc_{\Sigma}$ over $S$ as the mirror of $\Xc^\vee_{\Sigma}$.
Then $\Gamma(\ush_{\Lambda_{\Sigma}}\otimes_k \Phi(\Ec))$
should be interpreted as $\mathcal{W}(\Xc^\vee_\Sigma)$ by
descent properties of wrapped Fukaya categories developed in \cite{ganatrapardonshende2}, and
thus one should expect that
$\coh(\Xc_\Sigma) \simeq \mathcal{W}(\Xc_\Sigma^\vee)$.
This relationship is summarized as follows.
\[
\begin{tikzcd}
X_\Sigma \arrow[r, hook] \arrow[drr, dashed, <->]
& \Xc_\Sigma \arrow[r, dashed, <->] \arrow[d]
& \Xc_\Sigma^\vee \arrow[d]
& \arrow[l, hook'] S^\vee \arrow[dll, dashed, <->] \\
& S & X_\Sigma^\vee 
\end{tikzcd}
\]
Therefore mirror symmetry in this case
exhibits a fiber-base duality.

As a sanity check, we conclude our discussion with an example.
In \cite{hanlon19}, a particular type of
monodromy action on the Fukaya-Seidel category
is studied and is shown to be mirror to
action of the Picard group on the mirror toric
variety.
Namely, consider a toric variety specified by
the fan data $\Xi$ on the lattice \footnote{Here we avoid using the letters $\Sigma$ and $N$, because, as we will see, the fan data here play a different role than the fan data in the previous discussions do.} $K$,
and let $A \subseteq K\simeq \Z^n$ be the primitive generators.
The mirror Landau-Ginzburg model
is given by the superpotential
\[
W_\Xi = \sum_{\alpha \in A} c_{\alpha}z^\alpha,
\]
where $c_\alpha \in \C^\times$.
\begin{theorem}[\cite{hanlon19}]\label{hanlon-theorem}
The monodromy action on the Fukaya-Seidel category given by varying $c_\alpha$
around $0$ for one loop corresponds to
$- \otimes \Oc(D_\alpha)$ on $\coh(X_\Xi)$.
\end{theorem}
Let $B \subseteq A$ and write $l = \lvert B \rvert$.
Note that by viewing $c_\beta$ for $\beta\in B$ as variables in the above formula,
$W_\Xi$ can be extended to a function on
$(\C^\times)^n \times (\C^\times)^B \simeq (\C^\times)^{n + l}$,
denoted by $\widetilde{W_\Xi} : (\C^\times)^{n + l} \rightarrow \C$.
One sees that the Liouville sector $\C^{n + l}\backslash \widetilde{W_\Xi}^{-1}(\infty)$
is fibered over $(\C^\times)^l$ with fiber $(\C^\times)^{n}\backslash W_\Xi^{-1}(\infty)$ and Theorem \ref{hanlon-theorem} says that the monodromy
action map $\pi_1(\C^\times) \rightarrow \pic(X_\Xi)$ is given by sending the $i$-th generator to $\Oc(D_{\alpha_i})$.

On the other hand, note that the mirror of $\C^\times \simeq T^*S^1$ is $\A^1 - 0$.
Let $\Sigma$ be the fan of $(\A^1 - 0)^l$, consisting of only the origin inside $\Z^l$,
and set
\[\Ec = (\Oc(-D_{\beta_1}), \dots, \Oc(-D_{\beta_l})).
\]
The principal toric fibration $\Xc_{\Sigma,\Ec} \rightarrow X_\Sigma$  constructed
from these data is nothing but
the fiber product
\[
E(\Oc(D_{\beta_1}))\backslash Z_{X_\Xi} \times_{X_\Xi} \cdots \times_{X_\Xi} E(\Oc(D_{\beta_l}))\backslash Z_{X_\Xi}
\]
of the total spaces of the line bundles $\Oc(D_{\beta})$ with zero sections removed.
Then we have $\Xc_{\Sigma,\Ec}$ being mirror
to $((\C^\times)^{k + l}, \widetilde{W_\Xi})$.

\subsection{Outline of the paper}
In Section 2, we recall the necessary background on $\oo$-categories and sheaf theory.
In Section 3, we introduce the subjects of study of this paper -- principal
toric fibrations and the associated local systems of coefficients,
among other relevant definitions necessary to state our theorem.
Section 4 is devoted to the proof of the theorem in the affine case
-- when the toric fiber is $\A^k\times (\A^1 - 0)^{n-k}$,
which in Section 5 will be further generalized
to the case when the fiber
is a quotient stack $[\A^k\times (\A^1 - 0)^{n-k} / H]$ for some finite abelian group $H$.
We then show in Section 6 that
the general case can be obtained
from the affine case by a descent
argument, \emph{i.e.} gluing from
the local pieces where the fibers
are affine.
In Section 7, we provide an alternative presentation of our theorem, \`a la Harder-Katzarkov,
when the toric fiber is $\mathbb{P}^n$.

\section*{Acknowledgements}
YH thanks Qingyuan Bai for detailed discussions on the technical side of this work, thanks Wenyuan Li, and Dogancan Karabas for discussions on sheaf theory, thanks Diogo Fonseca for discussions on $\infty$-categories, thanks Qingyuan Bai, Bohan Fang, and Mingyuan Hu for various discussions through the years, thanks Zhenyi Chen and Ben Zhou for their interest in this work before the preparation of this manuscript. YH thanks Bohan Fang for hosting him as a visitor at Peking University the time the project was initiated.
YH thanks Eric Zaslow
who carefully read an early draft 
and provided many valuable suggestions.
Finally, YH wants to express his gratitude to Gus Schrader,
for reading multiple drafts and providing innumerable
suggestions and feedback,
without whose constant encouragement and support
this paper would not have reached its present form.

This work stems from the second author's Ph.D. thesis. PS sincerely thanks his advisor, Eric Zaslow, for his patience and guidance. Also, PS thanks Nicol\`{o} Sibilla and Hiro Lee Tanaka for answering many questions regarding quiver representations with coefficients in categories.

%% file: sections/notations.tex
\section{Notations and generalities}

\input{sections/notations/categories}
\input{sections/notations/sheaves}
\section{Principal toric fibrations}
\input{sections/notations/fibration}
\input{sections/notations/coefficient-system}

\subsection{Main theorem}
With these technical ingredients in place,
we may now state our main theorem.
\begin{theorem}\label{maintheorem}
There is an equivalence of $k$-linear stable $\infty$-categories,
\[
    \ind\coh \mathcal{X}_{\hat\Sigma,\beta,\mathcal{E}} \simeq \Gamma\left(\Lambda_{\hat{\Sigma},\beta}, \ush_{\Lambda_{\hat\Sigma}}^\diamond \otimes_k \Phi(\mathcal{E})\right).
\]
\end{theorem}
\begin{remark}
Since the sheaf is supported on $\Lambda_{\hat\Sigma, \beta}$, there is no difference between taking global sections over $\Lambda_{\hat\Sigma,\beta}$, over $T^*M_T$, and over $M_T$.
\end{remark}

%% file: sections/notations/categories.tex
\subsection{Categories}
\subsubsection{dg categories}
Let $k$ be an algebraically closed field of characteristic $0$.
Throughout this paper, we work with (the underlying $\oo$-category of) the category $dgCat_k$ of small dg categories
over $k$ with the Morita model structure.
In other words, we consider such categories only up to Morita equivalence. In particular, we only
consider idempotent-complete small dg categories.\footnote{This is possible because a small dg-category is always Morita equivalent to its idempotent completion.}
It turns out this is equivalent to considering (compactly-generated presentable) $k$-linear
stable $\infty$-categories because of the following result.
\begin{theorem}[{\cite[Corollary 5.7]{cohn2013differential}}]
There is an equivalence of $\infty$-categories
\[
N(\dgCat)[W^{-1}] \simeq \Mod_{Hk-mod}((Pr_{st,\omega}^L)^\otimes).
\]
Here the left-hand side is the $\infty$-localization of the category of
small $k$-dg-categories at Morita equivalences
(i.e.,  the underlying $\oo$-category of said model category) while the righ-hand side is
the $\infty$-category of compactly-generated presentable $k$-linear
stable $\infty$-categories with functors that preserve compact objects and colimits.
\end{theorem}
\begin{remark}
As noted in \cite[10.3]{gr1}, this is also
the $\infty$-category $\DGCat$ of (large) cocomplete $k$-dg-categories with continuous functors preserving compact objects.
Given an exact functor between idempotent complete small dg-categories, taking ind-completion gives a continuous functor between cocomplete dg-categories that also preserves compact objects.
Moreover the functor sending a small dg-category to its ind-completion is in fact
an equivalence of $\oo$-categories
\[
\mathbb{\ind}: \dgCat^{\text{ex}} \xrightarrow{\sim} \DGCat^{\text{cont}}.
\]
\end{remark}
\begin{remark}
Here $\Dc(k) := Hk-mod$ is the derived stable $\infty$-category of $k$-modules,
which is cocomplete and carries a symmetric monoidal structure
making it a commutative algebra in the sense of
\cite{lurie2017higher}.
\end{remark}
The theory of $\infty$-categories provides a natural
habitat for us to consider things in a homotopy coherent way,
which turns out to be extremely useful
when it comes to the study of functorial behavior of assignments of dg-categories.
The canonical references are \cite{lurie2009higher}\cite{lurie2017higher}.
For a version adapted to our usage,
we follow the treatments in \cite[Chapter 1]{gr1}.

It is often useful to switch between the two viewpoints: dg-categories and $k$-linear stable $\infty$-categories,
without explicitly indicating so, and we will do so from now on.
For convenience, we also write $\Cat_k$ instead of $\Mod_{Hk-mod}((Pr_{st,\omega}^L)^\otimes)$ for the $\infty$-category
of small $k$-dg categories / cocomplete $k$-dg-categories with continuous functors preserving compact objects.

\subsubsection{Limits and colimits}
When we say limits/colimits of diagrams of dg categories, i.e. functors
$p : K \rightarrow \Cat_k$,
we always mean limits/colimits in the $\infty$-categorical sense.
This is equivalent to considering homotopy limits/colimits in the model category sense,
when $K$ is actually (the nerve of) a $1$-category.

The following lemma says that under certain conditions, limits can be computed by colimit and vice versa, which will be useful later.
\begin{lemma}[{\cite[Corollay 5.5.3.4]{lurie2009higher}\cite[Proposition 1.7.5]{drinfeld2015compact}\cite[I, Chapter 1, 8.4.3]{gr1}}]
\label{limit-colimit}
Let $I$ be a small category, $\Psi:I\rightarrow \Cat_{k}$ be a functor, and write $C_i=\Psi(i)$. Suppose that for any $i\rightarrow j$, $\Psi(i\rightarrow j)$ admits a left adjoint $\phi_{ij}$.
Then there is a functor $\Phi: I^\op\rightarrow \Cat_k$ with $\Phi(i)=C_i$, $\Phi(i\rightarrow j)=\phi_{ji}$.
Moreover, there is a canonical equivalence
\[
\colim\Psi = \colim_{i\in I}C_i \simeq \lim \Phi = \lim_{i\in I^\op} C_i.
\]
\end{lemma}

\subsubsection{Tensor products of dg categories}\label{subsec:tensor-products}
The category $\Cat_k$ has a symmetric monoidal structure defined by Lurie\cite[4.5.2]{lurie2017higher}.
Here we follow the exposition in \cite[1.10.3]{gr1}.
By definition, the tensor product of two dg categories
$\mathcal{C}$ and $\mathcal{D}$ is $\mathcal{C}\otimes_k\mathcal{D} := \mathcal{C}\otimes_{Hk-mod}\mathcal{D}$
the relative tensor product where $\Cc$ and $\Dc$ are viewed as $Hk-mod$ modules.
In particular, the category of sheaves valued
in $\Cat_k$ inherits a tensor product,
also denoted by $-\otimes_k-$.
The following lemma will be useful later.
\begin{lemma}\label{mod-tensor}
Let $\Gamma$ be a 1-category, $\Cc$ be a cocomplete dg category. There is a canonical equivalence \[
\map(\Gamma, \Dc(k)) \otimes_k \Cc \simeq \map(\Gamma, \Cc).\]
\end{lemma}
\begin{proof}
In \cite[1.8.5.7]{gr1}, take $\mathcal{A}=k\Gamma$ the $k$-linearization of $\Gamma$ and $\mathbf{M} = \Cc$.
\end{proof}

%% file: sections/notations/sheaves.tex
\subsection{The Kashiwara-Schapira stack}
Here we give a brief account of the theory of constructible and microlocal sheaves adapted to our setting.
In this paper, the sheaves we consider are also known as constructible sheaves of unbounded chain complexes of $k$-modules.

\begin{definition}
A \emph{presheaf} of $k$-modules on $M$ is a functor \[\mathcal{F}\colon\open(M)^\op\rightarrow \Dc(k),\]
and we write $\psh(M;\Dc(k))$ for the category
of presheaves of $k$-modules.
\end{definition}
\begin{definition}
By a \emph{sheaf} of $k$-modules, we mean a
presheaf valued in $\Dc(k)$ satisfying (\v{C}ech) descent.
\end{definition}
When there is no confusion, we omit $\Dc(k)$
and simply write $\psh(M)$ and $\sh(M)$.
We may also write $\sh^\diamond(M)$ to
emphasize that we're working with sheaves of unbounded complexes.
\begin{remark}
This definition agrees with the classical approach of considering
the derived category of sheaves of abelian groups.
\end{remark}
\begin{remark}
The inclusion $\sh(M;\Dc(k)) \rightarrow \psh(M;\Dc(k))$ admits a left adjoint, the sheafification functor, denoted by $a:\psh(M;\Dc(k)) \rightarrow \sh(M;\Dc(k))$, exhibiting $\sh(M;\Dc(k))$ as an accessible localization of $\psh(M;\Dc(k))$.
\end{remark}
Most of the foundational results can be found in \cite{kashiwaraschapira} and some of the results adapted in the $\infty$-categorical
setting can be found in \cite{robalo2018lemma}.
\begin{remark}
In \cite{robalo2018lemma}, the authors consider hypercomplete sheaves.
However, this requirement is redundant in our setting, because every sheaf of $k$-modules is hypercomplete on a manifold.
This follows from the fact that a sheaf of $k$-modules is hypercomplete if and only if the underlying sheaf of spectra is hypercomplete,
if and only if the sheaf of spaces obtained by taking $\Omega^\infty$ is hypercomplete. However, sheaves of spaces on a manifold are always hypercomplete,
see \cite[7.2]{lurie2009higher}.
\end{remark}
\begin{definition}
Let $\mathcal{S} = \coprod S_\alpha$ be a stratification of $M$.
An $\mathcal{S}$-constructible sheaf of $k$-modules is a sheaf $\mathcal{F}$ such that $\mathcal{F}|_{S_\alpha}$ is a local system for any $\alpha$.
\end{definition}
Let $\mathcal{F}$ be an $\mathcal{S}-$constructible sheaf.
For any codirection $p = (x, \xi) \in T^*M$,
there is the \emph{microlocal stalk} $\ustk_p(\mathcal{F})$ of $\mathcal{F}$ at $(x, \xi)$.
To describe the functor, we choose a (germ of)
Morse function $f$ such that $df|_x = \xi$ and a small
neighborhood $U$ of $x$.
The microlocal stalk is the cofiber of
\[
\mathcal{F}(U\cap f^{-1}(\epsilon))
\rightarrow
\mathcal{F}(U\cap f^{-1}(-\epsilon)),
\]
for small $\epsilon > 0$.
Roughly speaking, the microlocal stalk at a codirection 
measures the failure of a sheaf to be constant
along said codirection.
\begin{definition}
The \emph{singular support} $\ssupp(\mathcal{F})$ of $\mathcal{F}$
is the closure of the set of codirections $p$ such that $\ustk_p(\mathcal{F})$ is non-zero.
\end{definition}
Denote by $\Lambda = \Lambda_\mathcal{S} := N^*\mathcal{S} = \coprod N^* S_\alpha$ the conormal of the stratification.
To each stratum $i_\alpha: S_\alpha \hookrightarrow M$, one can associate a sheaf $k_\alpha:=i_{\alpha*}k$.
We have the following.

\begin{lemma}
$\ssupp(k_\alpha) \subseteq \Lambda$. \qed
\end{lemma}

Viewing $\mathcal{S}$ as a poset where there
is a map $\alpha \rightarrow \beta$ iff
$S_\alpha \subseteq \overline{S_\beta}$,
the restriction maps together with the assignments $\alpha \mapsto k_\alpha$ give rise to a functor $\mathcal{S}^\op \rightarrow \sh_\Lambda(M)$.
The following is well-known, see for example \cite[3.1]{ganatrapardonshende2}.
\begin{theorem}
If all strata in $\mathcal{S}$ are contractible, then $\sh_\Lambda(M) \simeq \fun(\mathcal{S}, \Dc(k))$.
The equivalence is given by the following left Kan extension.
\[
\begin{tikzcd}
\mathcal{S}^\op \arrow[r] \arrow[d, "\yo"] & \sh(M) \\
\fun(\mathcal{S}, \Dc(k)) \arrow[dashed]{ur}{\simeq}[swap]{\text{LKan}} & \\
\end{tikzcd}
\]
\end{theorem}
In other words, the sheaf category $\sh_\Lambda(M)$ is compactly generated
by the corepresentables $k_\alpha$, the endomorphism algebra of which is $k[\mathcal{S}^\op]$.

\begin{theorem}\cite[Theorem 4.13]{ganatrapardonshende3}
Suppose $X \subset T^*M$ is closed conical, and $\Lambda \subset T^*M\backslash X$ is closed conical subanalytic isotropic.
The inclusion $\sh_X(M) \rightarrow \sh_{X \cup \Lambda}(M)$ is the full subcategory of $S$-right-orthogonal objects
and hence admits a left adjoint exhibiting the former as a localization.
Here $S$ is the collection of corepresentatives of microlocal stalk functors at the Lagrangian points of $\Lambda$.
\end{theorem}

One can define a sheaf of categories $\ush$, usually referred to as
\emph{microlocal sheaves} or the Kashiwara-Schapira stack,
as the sheafification of the following presheaf on $T^*M$
with the conic topology:
\[
\ush^{pre}: U \mapsto \sh(M)/\sh_{T^*M\backslash U}(M).
\]
This definition goes back to \cite[Section 6]{kashiwaraschapira} and plays
an important role in the study of the so-called microlocal hom space (see \cite[Definition 4.1.1]{kashiwaraschapira}) $\uhom$ between sheaves.
We have so far largely omitted this aspect of the theory
but want to remind the reader of the following theorem
as a motivation for the definition.
\begin{theorem}\cite[Theorem 6.1.2]{kashiwaraschapira}\cite[Corollary 5.5]{guillermou2012quantization}
Let $\Fc$ and $\mathcal{G}$ be sheaves on $M$ and denote by the same names their
images in $\ush(U)$.
There is a canonical equivalence
\[
\hom_{\ush(U)}(\Fc, \mathcal{G}) \simeq \uhom(\Fc, \mathcal{G}).
\]
\end{theorem}
For $\mathcal{F}\in \ush(U)$, there is a well-defined
notion of singular support which allows us to
consider microlocal sheaves ``supported" on $\Lambda$.
\begin{definition}
The \emph{Kashiwara-Schapira} stack $\ush_\Lambda$ associated to $\Lambda$ is the subsheaf of $\ush$
given by
\[
\ush_\Lambda: U \mapsto \{\mathcal{F} \in \ush(U): \ssupp(\mathcal{F})\subseteq \Lambda\}.
\]

\end{definition}

%% file: sections/notations/fibration.tex
\subsection{Stacky principal toric fibrations}\label{subsection:fibration}
In this paper, we consider a special case of fibrations $\mathcal{X}\rightarrow Y$ with toric fibers,
which we call \emph{stacky principal toric fibrations}.
A stacky principal toric fibration is a relative version of a toric stack
\cite{borisov2005orbifold, iwanari2009category, fantechi2010smooth, tyomkin2012}
(more precisely, a toric stack in the sense of \cite{tyomkin2012}.)

We recall the notion of a stacky fan
used to define a toric stack.
\begin{definition}
A \emph{stacky fan} \cite{borisov2005orbifold} denoted by $(\hat\Sigma,\beta)$ is the following data
\begin{itemize}
    \item a map $\beta \colon L \rightarrow N$ between finite rank free abelian groups with finite cokernel,
    \item a fan $\Sigma$ inside $M_\R:=N_\R^\vee$ consisting of rational strictly convex cones,
    \item a fan $\hat{\Sigma}$ inside $L_\R^\vee$
    consisting of rational strictly convex cones.
\end{itemize}
We assume, in addition, as in  \cite{kuwagaki20},
that
\begin{notation}\label{combinatorial-equivalence}
$\beta_{\R}$ induces an isomorphism between two fans $\hat{\Sigma}$ and $\Sigma$.
\end{notation}
\end{definition}
Now $\beta$ induces a surjective morphism between (algebraic) tori
\begin{equation}
    \beta_{ k^\times}: L_{ k^\times} \rightarrow N_{ k^\times}.
\end{equation}
We set $G_{\beta}:=\ker{\beta_{ k^\times}}$.
The information of a stacky fan
allows one to define a toric stack
\[
X_{\Sigma, \beta} = \left[
\colim_{\hat\sigma \in \hat{\Sigma}} k[\sigma^\vee \cap L^\vee]/G_\beta
\right],
\]
where $\hat\sigma^\vee:=\{v\in L^\vee_\R: v(x) \geq 0, \forall x \in \sigma \}$ is the dual cone.

In order to generalize this construction
to the relative case, one needs to specify
how the fibers vary over the the base.
The case we want to consider is when
the monodromy can be encoded
by the following data.
\begin{definition}
An \emph{Ikari} \footnote{Meaning ``anchor(s)".} is a morphism of abelian groups $\Ec: L^\vee \rightarrow \operatorname{Pic}(Y)$.
\end{definition}
\begin{remark}\label{ikari-altenative-defn}
An Ikari is non-canonically equivalent to the choice of $n=\operatorname{rank} L$ line bundles (order considered) on $Y$.
\end{remark}
Given any sub-monoid $S$ of $L^\vee$, one has an $\Oc_Y$-algebra
\begin{equation}
    \Ec[S]:=\oplus_{s\in S}\Ec(s),
\end{equation}
where the algebra structure is given by the obvious morphism $\Ec(s_1)\otimes_{\Oc_Y}\Ec(s_2) \xrightarrow {\simeq} \Ec(s_1+s_2)$.
\begin{example}
Let $S = \Z_{\geq 0}^n \subset L^\vee = \Z^n$, $\Ec = (\mathcal{L}_1, \dots, \mathcal{L}_n)$. Then $\Ec[S] \simeq \Sym_{\Oc_Y}(\mathcal{L}_1\oplus \cdots\oplus\mathcal{L}_n)$, and $\Spec_{\Oc_Y}(\Ec[S])=E(\Lc_1^\vee \oplus \cdots \oplus \Lc_n^\vee)$, the total space of $\Lc_1^\vee \oplus \cdots \oplus \Lc_n^\vee$.

\end{example}
We then define
\begin{equation}
    X_{\hat{\Sigma},\Ec}:=  \colim_{\hat{\sigma}\in \hat\Sigma}\Spec_{\Oc_Y}\Ec[\hat\sigma^\vee  \cap L^\vee].
\end{equation}
Note that the algebra $\Ec[S]$ carries a canonical $G_\beta$-action described as follows.
\begin{lemma}
$G_\beta$ is canonically isomorphic to $\spec k[ \coker(\beta^\vee)]$.
\end{lemma}
\begin{proof}
Note that
\[
    \beta_{ k^\times}: L_{ k^\times}= \spec k[L^\vee]\rightarrow N_{ k^\times} = \spec k[N^\vee].
\]
$\ker\beta_{ k^\times}$ is then given by the fiber at identity, which is
\[
\spec k[L^\vee/N^\vee]=\spec k[ \coker(\beta^\vee)].
\]
\end{proof}
\begin{cor}
$\Ec[S]$ admits a canonical $G_\beta$-action given by
\[\Ec[S] \rightarrow \Ec[S]\otimes  k[ \coker(\beta^\vee)],
\]
where the map on the second factor is induced by $S\hookrightarrow L^\vee \xrightarrow{\pi}  \coker(\beta^\vee
)$.
\end{cor}
\begin{remark}
Alternatively, one can view $ \coker(\beta^\vee)$ as the character lattice of $G_\beta$.
Then the weight of the $G_\beta$-action
is given by the projection $S\rightarrow  \coker(\beta^\vee)$.
\end{remark}
\begin{definition}
The \emph{stacky principal toric fibration} associated to the stacky fan $(\hat\Sigma, \beta)$
and the Ikari $\Ec$ is the quotient stack
\begin{equation}
\Xc_{\hat{\Sigma},\beta,\Ec}:=
\left[
    X_{\hat{\Sigma},\Ec} / G_\beta
\right]
= \left[  \colim_{\hat{\sigma}\in \hat\Sigma}\Spec_{\Oc_Y}\Ec[\hat\sigma^\vee  \cap L^\vee] / G_\beta
\right].
\end{equation}
\end{definition}

%% file: sections/notations/coefficient-system.tex
\subsection{The coefficient system \texorpdfstring{$\Phi(\Ec)$}{Phi(E)}}\label{subsection:local-system}

Associated to an Ikari $\Ec$,
we introduce a local system $\Phi(\Ec)$ of ($k$-linear stable $\infty$-)categories on the topological torus $M_T=M\otimes_\Z (\R/\Z)$. Recall the following fact.
\begin{lemma}
For any (topological) torus $T$, there is a canonical equivalence $T\simeq B\pi_1(T)$.
\qed
\end{lemma}
Applying this equivalence to $M_T$, we get $M_T\simeq B\pi_1(M_T) \simeq BM$. Recall that we have maps
\[
M \xrightarrow{\beta^\vee} L^\vee \xrightarrow{-\Ec} \pic(Y),
\]
inducing maps on classifying spaces
\[BM\rightarrow BL^\vee \rightarrow B\pic(Y).\]
On the other hand, we consider the action of Picard group $\pic(Y)$ on the category of sheaves. Consider the groupoid $B\pic(Y)$ and the following assignment\[
\ast \mapsto \ind\coh (Y),
\]
\[
\ast \xrightarrow{\alpha} \ast \mapsto \ind\coh(Y)\xrightarrow{\otimes \mathcal{L}_\alpha} \ind\coh(Y),
\]
where $\mathcal{L}_\alpha$ is the line bundle in $\pic(Y)$ corresponding to the morphism $\alpha$ in $B\pic(Y)$.
\begin{prop}
This assignment underlies a functor
\[B\pic(Y) \rightarrow \Cat_k.\]
\end{prop}
\begin{proof}
In the same way, we have a functor of 1-categories \[B\pic(Y) \rightarrow dgCat_k,\] which sends the unique object to the dg category of coherent sheaves $\coh(Y)$.
Note that by restricting to the full-subcategory whose objects consist of exactly one for each isomorphism class, one can make everything strictly commutative and thus avoid any 2-categorical issues.
This functor is therefore well-defined.
Now the composite \[B\pic(Y) \rightarrow N({dgCat_k}) \rightarrow N(dgCat_k)[W^{-1}] \simeq \Cat_k\] is the desired functor.
\end{proof}
\begin{definition}\label{def:local-system}
The coefficient system $\Phi(\Ec)$ associated to $\Ec$
is given by the composite functor
\[M_T\simeq BM \rightarrow BL^\vee \xrightarrow{B(-\Ec)} B\pic(Y) \rightarrow \Cat_k.\]
\end{definition}

\begin{remark}
Informally, this is the local system given by the representation of the fundamental group $\pi_1(M_T)\simeq M \xrightarrow{\beta^\vee} L^\vee \xrightarrow{-\Ec} \operatorname{Pic}(Y) \xrightarrow{\otimes} \operatorname{Aut}(\ind\coh(Y))$.
\end{remark}

%% file: sections/affine.tex
%!TEX root = ../main.tex
\section{Affine case}

Our first goal is to prove the main theorem for non-stacky affine toric fibrations.
Let $\Sigma=\Sigma(\sigma)$ be a \emph{regular} fan consisting of faces of a single rational strictly convex cone $\sigma$ in $N$.
We are going to show the following.

\begin{theorem}{\label{affine}}
There is an equivalence of $k$-linear stable $\infty$-categories,
\[
    \ind\coh \mathcal{X}_{\Sigma(\sigma),\mathcal{E}} \simeq \Gamma\left(\Lambda_{\Sigma(\sigma)}, \ush_{\Lambda_{\Sigma(\sigma)}}^\diamond \otimes_k \Phi(\mathcal{E})\right).
\]
\end{theorem}

\subsection{Description of constructible sheaves}
We start with a local combinatorial description of the sheaf $\ush^\diamond_{\Lambda_\Sigma}$.

First note that the stalk of $\ush^\diamond_{\Lambda_{\Sigma}}$ at 0 is just $\sh^\diamond_{\Lambda_{\Sigma}}(M_{\mathbb{R}})$.
By regularity of $\Sigma$, the generators of $\sigma$ can be completed to (or are) a $\Z$-basis of $N$.
Note that the construction of FLTZ skeleton does not depend on the choice of a basis.
In particular, we may assume $N=\Z^n$ and $\Sigma$ is generated by $e_0,\dots,e_{k - 1}$.
Then locally near 0, the skeleton is given by $(\top)^k\times (-)^{n-k}$,
where $\top$ denotes $ \R_{\leq 0} \times \R \subseteq T^*\R$ and $-$ denotes the zero section $0\times \R \subseteq T^*\R$.
We can therefore use a K\"unneth argument to get the result.
\footnote{This argument is due to Vivek Shende, and we thank Wenyuan Li for communicating the idea and explaining the details.}

However, we opt to calculate the sheaf
of categories $\ush_{\Lambda_\Sigma}$
directly by a close inspection of the skeleton for the purpose of better illustration.
Observe that for each generator $\tau \in \Sigma(1)$, $\tau^\perp$ stratifies the space $M_\R$ into three parts, namely, left, center (which is $\tau^\perp$), and right, denoted by $\tau_l, \tau_c, \tau_r$, in which we orient in such a way that the covector $\tau\in M^\vee_\R$ points from left to right.
$M_\R$ then has an induced stratification, 
\[
\mathcal{S}=\coprod_{\phi\in\{l, c, r\}^{\Sigma(1)}}\bigcap_{\tau\in\Sigma(1)}\tau_{\phi(\tau)}.
\]
We hereafter identify strata of $\mathcal{S}$ with $\{l, c, r\}^{\Sigma(1)}$ when there is no confusion.
It is clear that $N^*\mathcal{S}$ contains $\Lambda_\Sigma$.
An $\mathcal{S}$-constructible sheaf is then represented by a functor $\mathcal{F}\in \map(\mathcal{S}, \Dc(k)) \simeq \sh_{\mathcal{S}}(M_{\mathbb{R}})$.

It remains to impose singular-support conditions on the functor category.
Observe that as a poset $\mathcal{S}\simeq(\athree)^{\Sigma(1)}$.
Let $\mu, \nu \in \mathcal{S}$, there is an arrow $\mu \rightarrow \nu$ if and only if there is a $W\subseteq \Sigma(1)$, such that
$\mu(\tau)=c, \nu(\tau)\in\{l,r\}$ for $\tau\in W$ and $\mu(\tau)=\nu(\tau)$ for $\tau\in \mathcal{S}-W$.
Suppose $\abs{W}=1$, $W=\{\tau_0\}$, if $\mu(\tau)=c, \nu(\tau)=l$, we claim that $\mathcal{F}(\mu)\simeq \mathcal{F}(\nu)$. To see this, choose any point $p\in \mu$.
Note that $\Lambda_{\Sigma}|_p = \bigcup_{\tau_0\supseteq \tau\in \Sigma}(-\tau)\times p$.
Thus the covector $(\tau_0, p) \notin \Lambda_{\Sigma}$ because $\sigma$ is strictly convex. Vanishing of the microlocal stalk at $(\tau_0, p)$ then implies the desired isomorphism.
Now for any $\mu\rightarrow\nu$,  if it factors to a sequence $\mu=\kappa_0\rightarrow \kappa_1\rightarrow \cdots \rightarrow \kappa_l=\nu$,
such that $\abs{W(\kappa_i, \kappa_{i + 1})}=1$, then $\mathcal{F}(\mu)\simeq \mathcal{F}(\nu)$.
Functors satisfying these relations are exactly those factoring through
\[
\mathcal{S} \simeq (\athree)^{\Sigma(1)}
\rightarrow (\atwo)^{\Sigma(1)}
\rightarrow \Dc(k).
\]
Conversely, it is straightforward to see that constructible sheaves satisfying such relations have singular support contained in $\Lambda_{\Sigma}$.
In conclusion, we have a pull-back diagram of categories,
\[
\begin{tikzcd}
\sh_{\Lambda_{\Sigma}}(M_{\R}) \arrow[r, "\sim"] \arrow[dr, phantom, "\lrcorner", very near start] \arrow[d]
& \map\left((\atwo)^{\Sigma(1)}, \Dc(k)\right) \arrow[d] \\
\sh_{\mathcal{S}}(M_{\R}) \arrow[r, "\sim"]
& \map\left((\athree)^{\Sigma(1)}, \Dc(k)\right).\\
\end{tikzcd}
\]
We therefore have proved,
\begin{lemma}{\label{local}}
The stalk of $\ush^\diamond_{\Lambda_\Sigma}$ at $0 \in T^n$ is given by $\map\left((\atwo)^{\Sigma(1)},\Dc(k)\right)$. Furthermore, pullback of representations along a point $\bullet \rightarrow (\atwo)^{\Sigma(1)}$ is exactly the restriction functor to its respective stratum.\qed

\end{lemma}

One can cover $S^1=\R/\Z$ by two open intervals $\mathcal{U} = \{U, V\}$, say $U=(-1/4, 1/4)$, $V=(1/8, 7/8)$, whose \v{C}ech nerve is given by the following diagram $I = C(\mathcal{U})$:
\[
\begin{tikzcd}[sep=small]
& U & \\
*_l\arrow[ur] \arrow[dr] & & *_r \arrow[ul] \arrow[dl] \\
& V & \\
\end{tikzcd}.
\]

The products of these intervals give rise to a cover of $T^n\simeq M_T$, of which the \v{C}ech nerve is the previous diagram to the power of $n$.
A vertex in the diagram is labeled by an element $\kappa \in \{U, V, l, r\}^{[n]}$, denoted by $O_{\kappa}=\prod_{i=0}^{n-1} \kappa(i) \subset T^n$.
By an argument completely analogous to the previous one, one can show that
\begin{lemma}
$\ush_{\Lambda_\Sigma}^\diamond(O_\kappa)\simeq \map\left((\atwo)^{\Sigma(1)\cap{\kappa^{-1}(U)}},\Dc(k)\right)$,
where we identify $\Sigma(1)$ as the subset $\{0, \dots, k - 1\}$ of $\{0, \dots, n - 1\}$. \qed
\end{lemma}
In other words, the functor $\ush_{\Lambda_{\Sigma}}^\diamond: I^n \rightarrow \cat_k$ is equivalent to the composition of $I^n \rightarrow \sset$ and $\map(-, \Dc(k)): \sset^\op \rightarrow \cat_k$, where the first $k$ components of the map is given by $k$-fold product of
\[
\begin{tikzcd}[sep=small]
& U & \\
*_l\arrow[ur] \arrow[dr] & & *_r \arrow[ul] \arrow[dl] \\
& V & \\
\end{tikzcd}
\rightarrow
\begin{tikzcd}[sep=small]
& \Delta^1 & \\
\Delta^0 \arrow[ur] \arrow[dr] & & \Delta^0 \arrow[ul] \arrow[dl] \\
& \Delta^0 & \\
\end{tikzcd},
\]
and the last $n-k$ components are collapsing to a point.
Therefore, $\Gamma(T^n, \ush^\diamond_{\Lambda_\Sigma}\otimes_k \Phi(\mathcal{E}))$ is given by the limit of the following diagram.
\[
\begin{tikzcd}[sep=small]
\map((\Delta^1)^{\Sigma(1)}, \ind\coh(Y)) \arrow[d] \arrow[r]
& \prod_i\map((\Delta^1)^{\Sigma(1)\cap [n] - i}, \ind\coh(Y))\\
\prod_i\map((\Delta^1)^{\Sigma(1)\cap [n] - i}, \ind\coh(Y))\arrow[ur, swap, "\prod_i \otimes \mathcal{L}_i"]
\end{tikzcd}
\]
% In order to compute the global section of a sheaf, it would be helpful to allow calculations using
% a certain sort of Grothendieck topology.
% \begin{lemma}
% A sheaf on a topological space is a sheaf with respect to the proper locally open embedding topology.
% \end{lemma}
\subsection{Coherent sheaves}
We collect some foundational results from \cite{gr1}, see \emph{op. cit.} for a comprehensive treatment.
\begin{lemma}[K\"unneth]
\cite[Chapter 4, Lemma 6.3.2]{gr1} If $X_1, X_2\in \operatorname{Sch_{aft}}$, then there is an equivalence of categories
\[
\boxtimes\colon \ind\coh (X_1) \otimes \ind\coh (X_2) \simeq \ind\coh(X_1\times X_2).
\]
\end{lemma}
\begin{lemma}[Zariski Descent]
\cite[Chapter 4, 4.2]{gr1} If $U=\coprod U_{\alpha}$ is a Zariski cover of $Y$, then there is an equivalence
\[
\ind\coh(Y) \simeq \lim_{I \in C(U)} \ind\coh (U_{I}).
\]
\end{lemma}
Let $U=\coprod U_{\alpha}\rightarrow Y$ be a Zariski cover trivializing the line bundles $\mathcal{L}_i$.
Write $\mathcal{X}|_{U_{\alpha}}\cong (\A^k\times \Gm^{n - k}) \times U_{\alpha} $.
By K\"unneth,
\begin{equation}
\ind\coh(\mathcal{X}|_{U_\alpha})\simeq \ind\coh(\A^{k}\times \Gm^{n - k}) \otimes \ind\coh(U_{\alpha}).
\end{equation}
% Note that there is a canonical identification
% $\pi_*\mathcal{O}_{\mathcal{X}_{\alpha}} = \sym_{\mathcal{O}_{U_{\alpha}}}(\oplus \mathcal{L}_i)^\vee$
% under which $\pi_*p_\alpha^*(x_i)$ is identified with the 
% Consider the following commutative diagram
% \[
% \begin{tikzcd}[sep=small]
% p_\alpha^* \mathcal{O}_{\A^k\times \Gm^{n-k}}
% \arrow[rr, "p_{\alpha}^*x_i"]
% \arrow[d, equal]
% & & p_\alpha^* \mathcal{O}_{\A^\times \Gm^{n - k}}
% \arrow[r, equal]
% & \mathcal{O}_{\mathcal{X_{\alpha}}} \arrow[d, "\pi^*s^\alpha_i"] 
% \\
% \mathcal{O}_{\mathcal{X}_\alpha} \arrow[r, equal] & \pi^* \mathcal{O}_{U_\alpha} \arrow[rr, "\pi^*s_i^\alpha"]  & & \pi^*\mathcal{L}_i|_{{U}_\alpha},
% \end{tikzcd}
% \]
Note that we have the following equivalence \footnote{We thank Qingyuan Bai for extensive discussions and essentially coming up with this proof.}
\begin{lemma}
\[
    \ind\coh(\A^k\times \Gm^{n - k}) \simeq \map(\mathfrak{B}\Z_{\geq 0}^k \times \mathfrak{B}\Z^{n - k}, \Dc(k)),
\] where $\mathfrak{B}\Z$ (resp. $\mathfrak{B}\Z_{\geq 0}$) is the
category with a single object
whose endomorphism space is given by
$\Z$ (resp. $\Z_{\geq 0}$).\footnote{
We use $\mathfrak{B}\Z_{\geq 0}$ instead of $B\Z_{\geq 0}$ to stress that it is
a category but \emph{not} (necessarily) an $\infty$-groupoid, as the latter
is usually used to denote its 
$\infty$-groupoid completion.
In particular, $\mathfrak{B}\Z_{\geq 0} \centernot\simeq \mathfrak{B}\Z \simeq B\Z \simeq B\Z_{\geq 0}$.
}
\end{lemma}
\begin{proof}
We assume for simplicity that $n=k=1$, as the general case follows from the exact same argument or K\"unneth.
We therefore need to identify $\ind\coh(\A^1)$ with $\map(\mathfrak{B}\Z_{\geq 0}, \Dc(k))$.
The former can be identified with $\Dc(k[x])$,
while the latter as a functor category is cocomplete as $\Dc(k)$ is.
Note that $\map(\mathfrak{B}\Z_{\geq 0}, \Dc(k)) \simeq \map_k (k\mathfrak{B}\Z_{\geq 0}, \Dc(k))$ has a compact generator denoted by $\mathbf{1}$
given by the (enriched) Yoneda embedding
of the point $\ast\in (\mathfrak{B}\Z_{\geq 0})^\op$.
One sees that $\operatorname{End}(\mathbf{1}) \simeq k[x]$ and the result follows from Schwede-Shipley recognition principle (see \cite[7.1.2.1]{lurie2017higher}).
\end{proof}
Therefore
\begin{equation}
    \ind\coh(\mathcal{X}|_{U_\alpha}) \simeq \map(\mathfrak{B}\Z_{\geq 0}^k \times \mathfrak{B}\Z^{n - k}, \ind\coh(U_\alpha))
\end{equation} by Lemma \ref{mod-tensor}.
Note that we can also write the latter in terms of the interval groupoid $I$ as
\begin{equation}
    \map((\Delta^1/\partial \Delta^1)^k\times (I/\partial I)^{n - k}, \ind\coh(U_\alpha)),
\end{equation}
or equivalently as the limit of the diagram
\[
\begin{tikzcd}[sep=small]
\map((\Delta^1)^k \times I^{n - k}, \ind\coh(U_\alpha)) \arrow[d] \arrow[r]
& \prod_i\map(K_i, \ind\coh(U_\alpha))\\
\prod_i\map(K_i, \ind\coh(U_\alpha))\arrow[ur, swap, "\operatorname{id}"],
\end{tikzcd}
\]
where $K_i$ is the diagram by replacing the $i$-th component of $(\Delta^1)^k \times I^{n - k}$ by a point,
and the vertical (resp. horizontal) map is given by pullback along $K_i\rightarrow (\Delta^1)^k \times I^{n - k}$ the inclusion of the left (resp. right) endpoint in the $i$-th coordinate.
Although this indentification is not natural in $U_{\alpha}$,
it can be made so by tensoring with a twist
\[
\phi \in \map((\Delta^1)^k\times I^{n - k}, \ind\coh(Y))
\] defined as follows.
Denote a vertex $\lambda \in (\Delta^1)^k \times I^{n - k}$ by a sequence $\{\lambda_i\}$, where $\lambda_i = 0, 1$.
Let $\phi(\lambda) := \otimes_i \mathcal{L}_i^{-\lambda_i}$,
and the morphism $\mathcal{O}\rightarrow \mathcal{L}_i^\vee$ is given by the trivializing section.
The following is the key observation of this section.

\begin{prop}\label{prop:naturality}
After tensoring with $\phi$, the identification $\ind\coh(\mathcal{X}|_{U_\alpha}) \simeq \map((\Delta^1/\partial \Delta^1)^k\times (I/\partial I)^{n - k}, \ind\coh((U_\alpha))$ is independent of the choice of trivialization.
\end{prop}
\begin{proof}
Without loss of generality, we may assume $n = k = 1$.
Write $\mathcal{X}|_{U_{\alpha}} \simeq U_{\alpha} \times \A^1$ and  $p \colon \mathcal{X} \rightarrow Y$. Recall that $\Xc$ is the total space of $\Lc^\vee$.
Given a sheaf $\mathcal{F}$ on $\mathcal{X}|_{U_{\alpha}}$, the corresponding diagram is given by $p_*\mathcal{F} \xrightarrow{p_*x} p_*\mathcal{F}$.
Tensoring with $\phi$, it becomes
\[
\begin{tikzcd}
p_*\mathcal{F} \arrow[r, "{p_*x\otimes s_{\alpha}}"]
\arrow[d, equal]
& p_*\mathcal{F} \otimes \mathcal{L^\vee}
\arrow[d, equal]
\\
p_*\mathcal{F} \arrow[r, "p_*s"]
& p_*(\mathcal{F}\otimes p^*\mathcal{L^\vee})
\\
\end{tikzcd},
\]
where $s$ is the tautological section of $p^*\mathcal{L^\vee}$,
independent of the choice of trivialization.
\end{proof}
Note that the definition of $\phi$ encodes
exactly the monodromy of the local system $\Phi(\Ec)$.
The above proposition can thus be rephrased as the following.
\begin{prop}
Let $\Phi_\alpha(\Ec)$ be the local system defined by the same monodromy as that of $\Phi(\Ec)$
but with stalk $\ind\coh(U_\alpha)$.
Then there is an equivalence functorial in $U_\alpha$:
\[
\ind\coh(\mathcal{X}|_{U_\alpha}) \simeq \Gamma( \ush^\diamond_{\Lambda_{\hat{\Sigma}}}\otimes_ k\Phi_{\alpha}(\mathcal{E})).
\]
\qed
\end{prop}
We therefore have the following diagram functorial in $U_\alpha$,
\[
\begin{tikzcd}[sep=small]
\map((\Delta^1)^k \times I^{n - k}, \ind\coh(U_\alpha)) \arrow[d] \arrow[r]
& \prod_i\map((K_i, \ind\coh(U_\alpha))\\
\prod_i\map(K_i, \ind\coh(U_\alpha))\arrow[ur, swap, "\prod_i \otimes \mathcal{L}_i"].
\end{tikzcd}
\]
Now, the functoriality of this equivalence allows us to prove the theorem by gluing together
the local descriptions $\ind\coh(U_\alpha)$.
Taking limit in $U_\alpha$ gives the following by Zariski descent.
\[
\begin{tikzcd}[sep=small]
\map((\Delta^1)^k \times I^{n - k}, \ind\coh(Y)) \arrow[d] \arrow[r]
& \prod_i\map((K_i, \ind\coh(Y))\\
\prod_i\map(K_i, \ind\coh(Y))\arrow[ur, swap, "\prod_i \otimes \mathcal{L}_i"].
\end{tikzcd}
\]
The theorem then follows. \qed
\subsection{An abstract description}
For the cautious readers, we give an abstract
treatment to Proposition \ref{prop:naturality},
where the functoriality is evident.
\footnote{We thank Qingyuan Bai for suggesting this treatment.}
This section can be skipped by those content with previous arguments.

For notational simplicity, we consider the
case $n = k = 1$, remarking that the general
case follows almost verbatim.

Let $\Ac:=\Sym_{\Oc_Y}(\Lc)$ be the relative
symmetric algebra,
which is a commutative algebra object
in $\ind\coh(Y)$. We can therefore consider
$ \Mod_{\Ac}(\ind\coh(Y)) $,
the category of $\Ac$-modules in $\ind\coh(Y)$.

Write $\Xc = \Spec_{\Oc_Y}(\Ac) \xrightarrow{p} Y$ (the total space of $\Lc^\vee$). Given a sheaf $\Fc \in \ind\coh \Xc$, $p_*\Fc$ is canonically an $\Ac$-module, with the structure map given by adjunction
\footnote{
In particular, because $\Ac \simeq \oplus_{i\geq 0} \Lc^i$, we have a map $\Lc \otimes p_*\Fc \rightarrow p_*\Fc$, equivalently, a map
$p_*\Fc \rightarrow p_*\Fc \otimes \Lc^\vee$.
This is the same map as in the proof of Proposition
\ref{prop:naturality}.
}
\[
\Ac \otimes_{\Oc_Y} p_*\Fc \simeq p_*p^*p_*\Fc \rightarrow p_*\Fc.
\]
This gives rise to a functor
\footnote{
The adjunction pair
\[
\adjunction{p^*}{\ind\coh (\Xc)}{\ind\coh(Y)}{p_*}
\]
gives a monad $T = p_*p^*$ in $\ind\coh(Y)$.
The algebra object $\Ac$ is exactly the monad $T$
applied to the monoidal unit $\Oc_Y$.
Hence $\Mod_{\Ac}(\ind\coh(Y)) \simeq \Mod_T(\ind\coh(Y))$.
See \cite[4.7.3]{lurie2017higher}
for the construction of $p_*:\ind\coh(Y)\rightarrow \Mod_T(\ind\coh(\Xc))$.
}
\[
p_*: \ind\coh(\Xc) \rightarrow \Mod_\Ac(\ind\coh(Y)).
\]
Given an $\Ac$-module $\Mc$,
we have the diagram induced by $\Lc \rightarrow \Ac$
\[
\Lc\otimes\Mc \rightarrow \Mc.
\]
Tensoring with $\Lc^\vee$ on both sides gives
\[
\Mc \rightarrow \Mc\otimes \Lc^\vee.
\]
This gives a functor from $\Mod_\Ac(\ind\coh(Y))$ to
\[
\lim\left(
\begin{tikzcd}[sep=small]
\map(\Delta^1, \ind\coh(Y)) \arrow[d] \arrow[r]
& \map(\Delta^0, \ind\coh(Y))\\
\map(\Delta^0, \ind\coh(Y))\arrow[ur, swap, "- \otimes \mathcal{L}"]
\end{tikzcd}
\right) \simeq \Gamma(\ush^\diamond_{\Lambda_\Sigma}\otimes_k \Phi(\Ec)).
\]

Its composition with $p_*$ gives the desired functor
\[
\ind\coh(\Xc) \rightarrow \Gamma(\ush^\diamond_{\Lambda_\Sigma}\otimes_k \Phi(\Ec)).
\]
Now Proposition \ref{prop:naturality}
is a computation that says restricting to $U_\alpha$,
it is an equivalence.
Since both sides satisfy Zariski descent,
it is also an equivalence globally. This proves Theorem \ref{affine}.

% Moreover, it actually comes from a Zariski
% sheaf of categories on $Y$
% \[
% \Mod_\Ac \colon U \mapsto \Mod_{\Ac|_U}(\ind\coh(U)).
% \]
% Write $\Xc = \Spec_{\Oc_Y}(\Ac)$ (i.e. the total space of $\Lc^\vee$). There is another Zariski sheaf
% \[
% \ind\coh_\Ac \colon U \mapsto \ind\coh{\Xc|_U}.
% \]

%% file: sections/affine-stacky.tex
\section{Affine stacky case}

Consider the smooth affine toric fan $(\hat{\Sigma}, \beta)$ given by faces of a single rational strictly convex cone $\hat\sigma$ inside $L$.

The FLTZ skeleton in this case is given by
\[
\Lambda_{\hat\Sigma} = 
\bigcup_{\tau\in\hat\Sigma}\bigcup_{\chi \in M_{\tau,\beta}/M}(-\tau)\times (\tau^\perp + \chi) \subset T^*M_{T}.
\]

\subsection{Absolute case}
It is illustrative to first assume the base is $\spec k$ and focus on the stacky phenomena at the point.

\begin{lemma}
Let $\check{H}_\beta$ be the character lattice of $H_\beta$. We have $\check{H_\beta}\simeq M_{\sigma,\beta}/M$.
\end{lemma}

\begin{lemma}
\cite[lemma 7.2]{kuwagaki20}
\[
\Xc_{\hat\Sigma(\hat\sigma), \beta}:=[\spec k[\sigma^\vee \cap L^\vee]/G_\beta]\simeq[\spec k[\sigma^\vee/\sigma^\perp\cap M_{\sigma,\beta}]/H_\beta].
\]
\end{lemma}

\begin{lemma}
Under the identification $\check{H_\beta}\simeq M_{\sigma,\beta}/M$, the weight of an element in $ k[\sigma^\vee/\sigma^\perp\cap M_{\sigma,\beta}]$ is given by the projection $\sigma^\vee/\sigma^\perp\cap M_{\sigma,\beta} \hookrightarrow M_{\sigma,\beta} \rightarrow M_{\sigma,\beta}/M$.
\end{lemma}

Given a commutative monoid $\Lambda$ and a morphism of monoids $\pi \colon \Lambda \rightarrow \check{H}$,
where $\check{H}$ is the character lattice of a finite abelian group $H$, we define a category $\Gamma_{\Lambda, H}$ as follows.
\begin{itemize}
\item The objects are characters $\chi \in \check{H}$
\item For any two characters $\chi, \chi'$, $\hom(\chi, \chi') = \pi^{-1}(\chi' - \chi)$.
\end{itemize}
\begin{lemma}\label{stacky-coh-calculation}
Let $\Lambda$ be defined as above. With the action of $H$ on $\spec k[\Lambda]$ specified by $\pi$, we have
\[
\ind\coh [\spec k[\Lambda]/H] \simeq \map(\Gamma_{\Lambda, H}, \Dc( k)).
\]
\end{lemma}
\begin{proof}
Note that $\{\mathcal{O}_{\chi}:\chi\in \check{H} \}$ generates $\ind\coh [\spec k[\Lambda]/H]$.
Therefore we have
\[
\ind\coh [\spec k[\Lambda]/H] \simeq \Mod_{\Ac}(\Dc(k)),\]
where \[\Ac := \rhom(\oplus_{\chi \in \check{H}} \mathcal{O}_{\chi}, \oplus_{\chi \in \check{H}} \mathcal{O}_{\chi})\]
is the endomorphism algebra of $1 := \oplus_{\chi \in \check{H}}\Oc_{\chi}$.
A direct calculation \footnote{This is a special case of \cite[Proposition 7.4]{kuwagaki20} restated.} shows that
\[
\rhom(\Oc_{\chi_1}, \Oc_{\chi_2}) = k[\pi^{-1}(\chi_2 - \chi_1)].
\]

Hence $\Ac$ is the same as the dg-algebra formed from the $k$-linearization of $\Gamma_{\Lambda, H}$,
and thus
\[
\fun(\Gamma_{\Lambda, H}, \Dc(k)) \simeq \Mod_{\Ac}(\Dc(k)) \simeq \ind\coh [\spec k[\Lambda]/H].
\]

\end{proof}

In particular, we take $\Lambda = \Lambda(\hat\Sigma(\hat\sigma)):=\sigma^\vee/\sigma^\perp\cap M_{\sigma,\beta}$ and $H = H_{\beta}$.
On the other hand, we have the following description of sheaf categories.
\begin{lemma}\label{stacky-constructible-calculation}
\[
    \sh_{\Lambda_{\hat\Sigma}}^\diamond(M_T) \simeq \map(\Gamma_{\Lambda_{\hat\Sigma}, H_\beta}, \Dc(k))
\]
\end{lemma}
\begin{proof}
This can be obtained by gluing the local description
of $\ush^\diamond_{\Lambda_{\hat\Sigma,\beta}}$ around each character $\chi \in M_{\sigma, \beta}/M \subset M_T$. Or rather, it can be deduced from the fact that $\sh^\diamond_{\Lambda_{\hat\Sigma,\beta}}(M_T) \simeq \ind\coh \Xc_{\hat\Sigma(\hat\sigma)}$ as shown in \cite{kuwagaki20}.
\end{proof}
\begin{figure}
    \centering
    \includegraphics[width=0.3\textwidth]{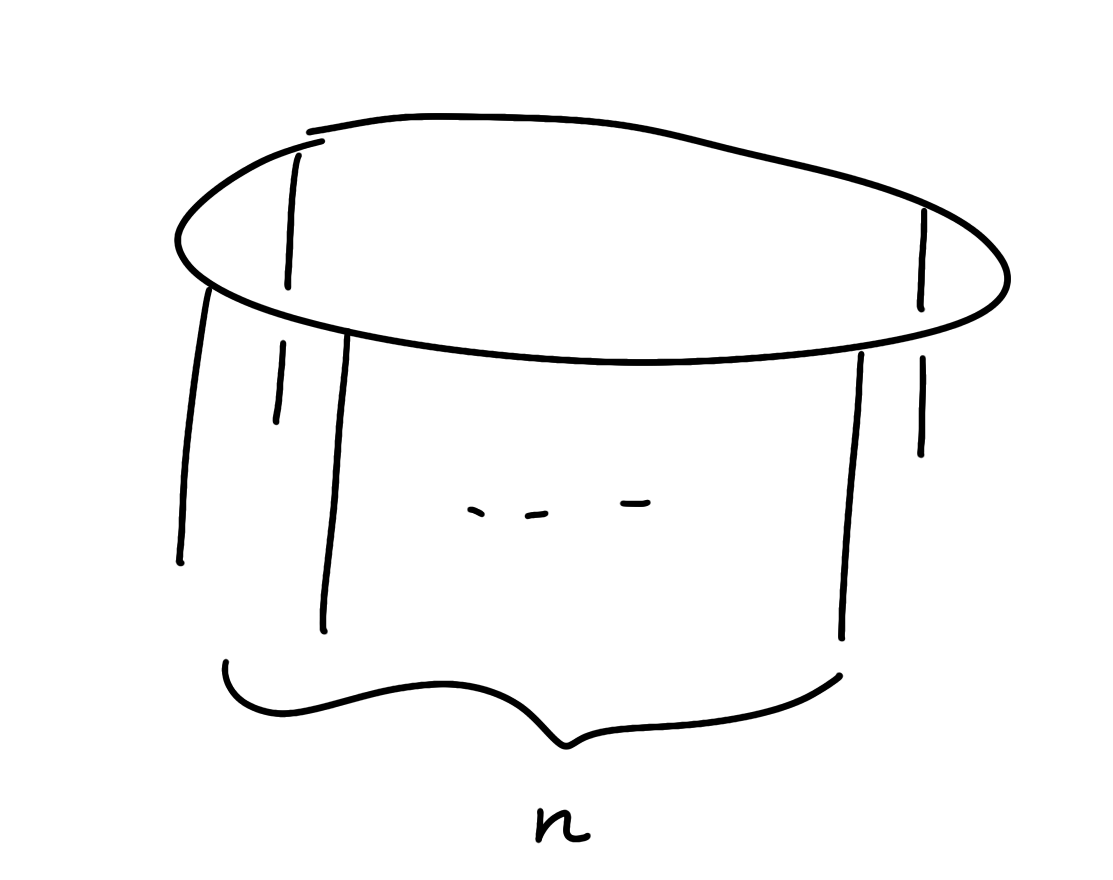}
    \caption{FLTZ skeleton for $[\A^1/\mu_n]$}
    \label{fig:skeleton-stacky-line}
\end{figure}
\begin{example}
If $\Xc = [\A^1/\mu_n]$, then $M_{\sigma, \beta} = \frac{1}{n}\Z$, $M = \Z$, and $H_\beta \simeq \check{H_\beta} \simeq \Z/n\Z$.
In this case, the monoid $\Lambda(\hat\Sigma) = \frac{1}{n}\Z_{\geq 0}$, and the map $\Lambda(\hat\Sigma) \rightarrow \check{H_\beta}$
is given by $\frac{i}{n} \mapsto i$.
Now the category $\Gamma_{\Lambda, H}$
has $n$ objects corresponding to characters of $\mu_n$.
The hom space between $i, j\in \Z/n\Z$ is given
by $(\frac{j - i}{n} + \Z)\cap \Z_{\geq 0}$.
Considering a cyclic quiver of $n$ vertices,
one sees that this hom space is the same
as the set of paths from $i$ to $j$,
and composition of morphisms in this category
corresponds to concatenation of paths.
Therefore $\Gamma_{\Lambda, H}$ is the
category freely generated by the cyclic quiver,
i.e., its path category.
A functor in $\fun(\Gamma_{\Lambda, H}, \Dc(k))$
is thus identified with a representation of the cyclic quiver.
A $\mu_n$-equivariant $k[x]$-module $M$
is then identified with the quiver representation,
where the vertex $\chi\in \Z/n\Z$ is mapped to $M_\chi$, the $\chi$-isotypic
component of $M$, and the arrow between $\chi$ and $\chi + 1$ is mapped to $M_{\chi}\xrightarrow{x}M_{\chi + 1}$.
On the other hand,
the category of constructible sheaves
is also identified with representations
of the same cyclic quiver, see figure $\ref{fig:skeleton-stacky-line}$.
\end{example}
\subsection{Relative case}

% \begin{lemma}
% Let $p: M \rightarrow M/G$ be a covering map, $\mathcal{F}$ a sheaf on $M/G$.
% Then $p^*\mathcal{F}$ has a natural $G$-equivariant structure, and $\Gamma(M, p^*\mathcal{F})^G = \Gamma(M/G, \mathcal{F})$.\qed
% \end{lemma}

% With this understanding, to prove the equivalence it suffices to compare two sheaves locally and identify the $M$-action.

Following previous notations, let $U=\coprod U_\alpha \rightarrow Y$ be a Zariski cover trivializing the line bundles.
Let $\Phi_{\alpha}(\mathcal{E})$ be the local system of categories defined by the same monodromy as that of $\Phi(\mathcal{E})$ but with stalk $\ind\coh(U_{\alpha})$, then $\Phi(\mathcal{E}) \simeq \lim_{\alpha} \Phi_{\alpha}(\mathcal{E})$.
\begin{remark}\label{rmk:trivialization}
The local system $\Phi_{\alpha}(\mathcal{E})$ is trivial, but not canonically so. The trivialization is given by the trivialization of the line bundles.
\end{remark}
\begin{prop}\label{local-local}
We have the following equivalence functorial in $U_\alpha$:
\[
\ind\coh(\mathcal{X}|_{U_\alpha}) \simeq \Gamma(M_T, \ush^\diamond_{\Lambda_{\hat{\Sigma},\beta}}\otimes_ k\Phi_{\alpha}(\mathcal{E})).
\]
\end{prop}
\begin{proof}
By K\"unneth and the absolute case, the LHS is equivalent to $\sh_{\Lambda_{\hat\Sigma}}(M_T)\otimes_ k \ind\coh(U_\alpha)\simeq \Gamma(M_T, \ush^\diamond_{\Lambda_{\hat{\Sigma}}}\otimes_ k \ind\coh(U_\alpha))\simeq \Gamma(M_T, \ush^\diamond_{\Lambda_{\hat{\Sigma}}} \otimes_ k \Phi_{\alpha}(\mathcal{E}))$.
The last equivalence is induced by the isomorphism in remark \ref{rmk:trivialization}.
By the proof of proposition \ref{prop:naturality}, this composed equivalence does not depend on the choice of trivialization and is functorial in $U_\alpha$.
\end{proof}
\begin{theorem}\label{thm:affine-stacky}
Let $\mathcal{X}_{\hat\Sigma(\sigma), \beta}$ be
the principal toric fibration defined by the stacky fan $(\hat\Sigma(\sigma), \beta)$ of a single rational strictly convex cone $\sigma$. Then we have
\[
\ind\coh(\mathcal{X}_{\hat{\Sigma}(\sigma),\beta}) \simeq \Gamma(M_T, \ush^\diamond_{\Lambda_{\hat\Sigma}}
\otimes_ k
\Phi(\mathcal{E})).
\]
\end{theorem}
\begin{proof}
Let $\{U_j\}_{j\in J}$ be a sufficiently fine finite cover of $M_T$ such that all intersections are (finite unions of) contractibles and that
\[
\Gamma(U_j, \ush^\diamond_{\Lambda_{\hat\Sigma,\beta}}\otimes_k \Phi_{\alpha}(\Ec)) \simeq \Gamma(U_j, \ush^\diamond_{\Lambda_{\hat\Sigma,\beta}}) \otimes_k \Gamma(U_j, \Phi_{\alpha}(\Ec)).
\]
Note that pullback $f^*$ of ind-coherent sheaves admits a right adjoint $f_*$ (cf. \cite[II, Chapter 4, 3.1.6]{gr1}), and restriction $i^*$ of microlocal sheaves to an open subset is cocontinuous.
Thus Lemma \ref{limit-colimit} applies and we can calculate as follows.
\begin{align*}
    \ind\coh(\Xc_{\hat{\Sigma}(\sigma),\beta}) & \simeq \lim_{\alpha\in C(A)}\ind\coh(\Xc|_{U_\alpha}) \\
    & \simeq \colim_{\alpha\in C(A)^\op}\ind\coh(\Xc|_{U_\alpha}) \\
    & \simeq \colim_{\alpha\in C(A)^\op} \Gamma(M_T, \ush^\diamond_{\Lambda_{\hat\Sigma,\beta}}\otimes_k \Phi_{\alpha})\\
    & \simeq \lim_{j\in C(J)}\colim_{\alpha\in C(A)^\op}\Gamma(U_j, \ush^\diamond_{\Lambda_{\hat\Sigma,\beta}})\otimes_k \Gamma(U_j, \Phi_{\alpha}(\Ec)) \\
    & \simeq \lim_{j\in C(J)} \Gamma(U_j, \ush^\diamond_{\Lambda_{\hat\Sigma,\beta}})\otimes_k\colim_{\alpha\in C(A)^\op} \Gamma(U_j, \Phi_{\alpha}(\Ec)) \\
    & \simeq \lim_{j\in C(J)} \Gamma(U_j, \ush^\diamond_{\Lambda_{\hat\Sigma,\beta}})\otimes_k\lim_{\alpha\in C(A)} \Gamma(U_j, \Phi_{\alpha}(\Ec)) \\
    & \simeq \lim_{j\in C(J)} \Gamma(U_j, \ush^\diamond_{\Lambda_{\hat\Sigma,\beta}})\otimes_k\Gamma(U_j, \Phi(\Ec)) \\
    & \simeq \lim_{j\in C(J)} \Gamma(U_j, \ush^\diamond_{\Lambda_{\hat\Sigma,\beta}}\otimes_k \Phi(\Ec)) \\
    & \simeq \Gamma(M_T, \ush^\diamond_{\Lambda_{\hat\Sigma,\beta}}\otimes_k\Phi(\Ec)).
\end{align*}
\end{proof}
% \textcolor{red}{
% For arbitrary toric stack $(\Sigma, \beta)$, if global section commutes with limit,
% $ Coh(Z) \cong \mathrm{lim}Coh(Z|_\alpha) \cong \mathrm{lim}\Gamma(\mu sh_{\Lambda_{\Sigma, \beta}}\otimes \Phi_\alpha) \cong \Gamma(\mu sh_{\Lambda_{\Sigma, \beta}} \otimes \mathrm{lim}\Phi_\alpha) \cong \Gamma(\mu sh_{\Lambda_{\Sigma, \beta}}\otimes \Phi)$ ?
% }
% \todo{yes}

%% file: sections/descent.tex
\section{Descent}
In this section, we are going to prove
our main theorem in the general case using
descent arguments.
Such a proof goes as follows: locally a
principal toric fibration is a Cartesian
product; for any Zariski open set on
which the fibration is trivial, one
can show that the theorem follows from
gluing along the fiber direction - i.e., from the local case where the toric fiber is affine;
next, we will show that we can further
glue along the base direction and use that to prove
the main theorem.
\subsection{Toric descent}
First, we are going to describe ``descent in the fiber direction", i.e., descent along the cover
given by the toric fan data.

Note that we can identify $\hat\Sigma$ with the Cech nerve $C(\hat\Sigma)$ of the set of its top-dimensional cones.
Let $\pmb{\sigma} = \{\sigma_0, \sigma_1, \dots, \sigma_k\} \subset \hat\Sigma_0$ be a collection of top-dimensional cones.
Write $\abs{\pmb\sigma} = \cap_i \sigma_i$, which is a codimension $k$ cone, and set
\[
\Lambda_{\pmb\sigma} := \Lambda_{\hat{\Sigma}(\abs{\pmb{\sigma}}), \beta}.
\]
We then have
\[
\Lambda_{\hat\Sigma, \beta} = \bigcup_{\pmb{\sigma}\in C(\Sigma)}\Lambda_{\abs{\pmb\sigma}}.
\]
Now given an inclusion of face $\sigma_1\subset \sigma_2$,
one can define a functor
\[
I^{\sigma_2\sigma_1}: \sh^\diamond_{\Lambda_{\hat\Sigma(\sigma_2),\beta}}(M_T)
\rightarrow
\sh^\diamond_{\Lambda_{\hat\Sigma(\sigma_1),\beta}}(M_T),
\]
given by
\[
(-) \mapsto (-) \ast \Theta(\sigma_1, 0),
\]
where $\Theta(\sigma, 0):=p_! k_{\interior(\sigma^\vee)}[n]$ is the ``costandard sheaf" associated to $\sigma$.

Recall the functor $\kappa$ as defined in \cite{fang2011categorification, treumann2010remarks, kuwagaki20}
given by
\[
\Theta'(\sigma, \chi) = j_{X_\sigma *}\Oc_{X_\sigma}(\chi) \mapsto \Theta(\sigma, \chi)= p_! k_{\interior{(\sigma^\vee + \chi)}}[n],
\]
we note the following.
\begin{lemma}
The identification established by combining lemmata \ref{stacky-coh-calculation}, \ref{stacky-constructible-calculation} is $\kappa$ composed with a degree shift by $-n$.
\end{lemma}
\begin{proof}
The sheaf
$\Theta'(\sigma, 0) = \Oc_{X_{\hat\Sigma(\hat\sigma), \beta}}$ is identified with $k[(\sigma^\vee/\sigma^\perp) \cap M_{\sigma,\beta}]$ as a $H_\beta$-equivariant module over itself. Viewed as an element in $\fun(\Gamma_{\Lambda_{\hat\Sigma}, H_\beta}, \Dc(k))$,
it is given by the formula
\[
\chi \mapsto \rhom(\Oc(\chi), \Oc) = k[(\sigma^\vee/\sigma^\perp) \cap M_{\sigma,\beta}]_\chi =  k[(\sigma^\vee/\sigma^\perp)\cap (\chi + M)],
\]
i.e., the $k$-module on the vertex $\chi$
is the $\chi$-isotypic component of $k[(\sigma^\vee/\sigma^\perp) \cap M_{\sigma,\beta}]$.

On the other hand, the stalk of the sheaf $p_!k_{\interior(\sigma^\vee)}$ at a character
$\chi \in \check{H_\beta} = M_{\sigma,\beta}/M\subseteq M_T$ is canonically
identified with the $k$-module
\[
k[(\sigma^\vee/\sigma^\perp)\cap(\chi + M)],
\]
by inspecting the fiber of $M_{\R}\rightarrow M_T$ at $\chi$.

We thus see that our identification sends $\Theta'(\sigma, 0)$ to $p_!k_{\interior(\sigma^\vee)}$.
Similarly, it sends $\Theta'(\sigma, \chi)$ to $p_!k_{\interior(\sigma^\vee + \chi)}$.
Therefore it differs from $\kappa$ by a degree shift
of $-n$.
\end{proof}
We have the following functoriality result.
\begin{theorem}[\cite{fang2011categorification, treumann2010remarks}{\cite[Corollary 12.8]{kuwagaki20}}]
The assignment
$\sigma\mapsto \sh^\diamond_{\Lambda_{\hat\Sigma(\sigma),\beta}}$ together with above defined $I$'s give rise to a functor $C(\Sigma)\rightarrow \Cat_k$.
Moreover, this assignment is compatible with
restrictions of ind-coherent sheaves.
In other words, for an inclusion of face $\tau\subset\sigma$, the following diagram commutes.
\[
\begin{tikzcd}
\sh^\diamond_{\hat\Sigma(\sigma),\beta}(M_T) \arrow[r, "\simeq"] \arrow[d, "{\ast\Theta(\tau,0)}"]
& \ind\coh(X_{\hat\Sigma(\sigma),\beta}) \arrow[d, "\iota_{\tau\sigma}^*"]\\
\sh^\diamond_{\hat\Sigma(\tau),\beta}(M_T) \arrow[r, "\simeq"]
& \ind\coh(X_{\hat\Sigma(\tau),\beta}).
\end{tikzcd}
\]
\end{theorem}
\begin{remark}
In \emph{op. cit.}, the horizontal arrows are $\kappa$, which differ
from the arrows shown above by a degree shift.
The only difference is that our equivalence
no longer intertwines tensor products and
convolution products for trivial reasons,
which is immaterial for our applications.
\end{remark}
We can now prove the following global version of Proposition \ref{local-local}.
\begin{prop}\label{global-local}
Assume that $\hat\Sigma$ is smooth. There is an equivalence
functorial in $U_\alpha$.
We have the following equivalence functorial in $U_\alpha$:
\[
\ind\coh(\mathcal{X}_{\hat\Sigma,\beta,\Ec}|_{U_\alpha}) \simeq \Gamma(M_T, \ush^\diamond_{\Lambda_{\hat{\Sigma},\beta}}\otimes_ k\Phi_{\alpha}(\mathcal{E})).
\]
\end{prop}
\begin{proof}
Similarly to Theorem \ref{thm:affine-stacky}, we have the following
sequence of equivalences.
\begin{align*}
\ind\coh(\Xc_{\hat\Sigma,\beta,\Ec}) & \simeq \lim_{\sigma\in C(\hat\Sigma)} \ind\coh(\Xc_{\hat\Sigma(\sigma), \beta, \Ec}) \\
& \simeq \colim_{\sigma\in C(\hat\Sigma)^\op} \ind\coh(\Xc_{\hat\Sigma(\sigma), \beta, \Ec})\\
& \simeq \colim_{\sigma\in C(\hat\Sigma)^\op} \ind\coh(X_{\hat\Sigma(\sigma), \beta})\otimes_k \ind\coh(U_\alpha) \\
& \simeq \left(\colim_{\sigma\in C(\hat\Sigma)^\op} \ind\coh(X_{\hat\Sigma(\sigma), \beta})\right)\otimes_k \ind\coh(U_\alpha) \\
& \simeq \left(\lim_{\sigma\in C(\hat\Sigma)} \ind\coh(X_{\hat\Sigma(\sigma), \beta})\right)\otimes_k \ind\coh(U_\alpha) \\
& \simeq \sh^\diamond_{\Lambda_{\hat\Sigma,\beta}}(M_T) \otimes \ind\coh(U_\alpha) \\
& \label{tensor-constant-sheaf}  \simeq \Gamma(M_T, \ush^\diamond_{\Lambda_{\hat\Sigma,\beta}}\otimes_k \ind\coh(U_\alpha)) \\
& \simeq \Gamma(M_T, \ush^\diamond_{\Lambda_{\hat\Sigma,\beta}}\otimes_k \Phi_\alpha(\Ec)).
\end{align*}
The last two rows are true because
one can take a sufficiently fine cover
and turn limits into colimits and back
as in the proof of Theorem \ref{thm:affine-stacky}.
As $\hat\Sigma$ is smooth,
Proposition \ref{local-local} applies to each local piece and we thus know the equivalence is functorial in $U_\alpha$.
\end{proof}
\subsection{Non-smooth case}
Assume now $\hat\Sigma(\hat\sigma)$ is affine but not smooth.
We choose a smooth refinement $\hat\Sigma(\hat\sigma) \subset \tilde{\Sigma}$
and denote the corresponding morphism on
toric stacks by $f:\Xc_{\tilde{\Sigma},\beta,\Ec} \rightarrow \Xc_{\hat\Sigma(\sigma),\beta,\Ec}$.
By abuse of notation, we also write the map between fibers by $f: X_{\tilde{\Sigma},\beta} \rightarrow X_{\hat\Sigma(\sigma),\beta}$.

Let $I:\sh^\diamond_{\Lambda_{\hat{\Sigma}(\hat\sigma)}}(M_T) \hookrightarrow \sh^\diamond_{\Lambda_{\tilde\Sigma}}(M_T)$ be the inclusion. 
\begin{theorem}[{\cite[Corollary 12.7]{kuwagaki20}}]\label{refinement-functoriality}
The following diagram is commutative.
\[
\begin{tikzcd}
\sh^\diamond_{\Lambda_{\hat{\Sigma}(\hat\sigma)}}(M_T) \arrow[r, hook, "I"] \isoarrow{d}
&  \sh^\diamond_{\Lambda_{\tilde\Sigma}}(M_T) \isoarrow{d} \\
\ind\coh(X_{\hat\Sigma(\sigma),\beta}) \arrow[r, "f^*"]
& \ind\coh(X_{\tilde\Sigma,\beta})\\
\end{tikzcd}
\]
\end{theorem}
We are now ready to prove:
\begin{prop}\label{singular-local-local}
Assume that $\hat\Sigma(\hat\sigma)$ is affine but not necessarily smooth. There is an equivalence
functorial in $U_\alpha$.
We have the following equivalence functorial in $U_\alpha$.
\[
\ind\coh(\mathcal{X}_{\hat\Sigma,\beta,\Ec}|_{U_\alpha}) \simeq \Gamma(M_T, \ush^\diamond_{\Lambda_{\hat{\Sigma},\beta}}\otimes_ k\Phi_{\alpha}(\mathcal{E})).
\]
\end{prop}
\begin{proof}
By Theorem \ref{refinement-functoriality},
we have the following commutative diagram,
\[
\begin{tikzcd}
\sh^\diamond_{\Lambda_{\hat{\Sigma}(\hat\sigma)}}(M_T) \otimes_k \ind\coh(U_\alpha)\arrow[r, "I\otimes 1"] \isoarrow{d}
&  \sh^\diamond_{\Lambda_{\tilde\Sigma}}(M_T)\otimes_k \ind\coh(U_\alpha) \isoarrow{d} \\
\ind\coh(X_{\hat\Sigma(\sigma),\beta})\otimes_k \ind\coh(U_\alpha) \arrow[r, "f^*\otimes 1"]
& \ind\coh(X_{\tilde\Sigma,\beta})\otimes_k \ind\coh(U_\alpha).\\
\end{tikzcd}
\]
Applying Lemma \ref{limit-colimit},
we can show that this is equivalent to
\[
\begin{tikzcd}
\Gamma(M_T, \ush^\diamond_{\Lambda_{\hat{\Sigma}(\hat\sigma)}} \otimes_k \ind\coh(U_\alpha))\arrow[r, "I\otimes 1"] \isoarrow{d}
&  \Gamma(M_T, \ush^\diamond_{\Lambda_{\tilde\Sigma}}(M_T)\otimes_k \ind\coh(U_\alpha)) \isoarrow{d} \\
\ind\coh(X_{\hat\Sigma(\sigma),\beta})\otimes_k \ind\coh(U_\alpha) \arrow[r, "f^*\otimes 1"]
& \ind\coh(X_{\tilde\Sigma,\beta})\otimes_k \ind\coh(U_\alpha).\\
\end{tikzcd}
\]
Note that the horizontal arrows are fully faithful by \cite[1.10.5.8]{gr1} because $I$ is.
Apply the equivalence $\Phi_\alpha(\Ec)\simeq \ind\coh(U_\alpha)$ given by the trivialization of the line bundles.
We have
\[
\begin{tikzcd}
\Gamma(M_T, \ush^\diamond_{\Lambda_{\hat{\Sigma}(\hat\sigma)}} \otimes_k \Phi_\alpha(\Ec))\arrow[r, "I\otimes 1"] \isoarrow{d}
&  \Gamma(M_T, \ush^\diamond_{\Lambda_{\tilde\Sigma}}\otimes_k \Phi_\alpha(\Ec)) \isoarrow{d} \\
\ind\coh(\Xc_{\hat\Sigma(\sigma),\beta,\Ec}) \arrow[r, "f^*"]
& \ind\coh(\Xc_{\tilde\Sigma,\beta,\Ec}).\\
\end{tikzcd}
\]
Now by Theorem \ref{global-local},
the right vertical arrow is an equivalence functorial in $U_\alpha$ and thus so is the left one by fully-faithfulness.
\end{proof}

\subsection{Proof of the main theorem}
Using Proposition \ref{singular-local-local} and
following the exact same line of Proposition \ref{global-local}
with the exception that our affine component can
now be singular, we prove that
\begin{theorem}
Let $(\hat\Sigma,\beta)$ be a stacky fan,
not necessarily smooth nor affine,
satisfying assumption \ref{combinatorial-equivalence}.
There is an equivalence of $k$-linear stable $\infty$-categories,
\[
    \ind\coh \mathcal{X}_{\hat\Sigma,\beta,\mathcal{E}} \simeq \Gamma\left(M_T, \ush_{\Lambda_{\hat\Sigma}}^\diamond \otimes_k \Phi(\mathcal{E})\right).
\]
\qed
\end{theorem}

%% file: sections/quiver.tex
\section{Quiver theoretic apporach for toric $\mathbb{P}^n$-bundles}\label{quiver}

\subsection{Quiver theoretic interpretation of the main theorem}
A version of CCC for toric bundles first appeared in the work of Harder-Katzarkov \cite{harder2019perverse}. The purpose of the remaining sections is to explain the relationship between their work and our main theorem. As a byproduct, we will find an alternative, quiver-theoretic proof of CCC for toric $\mathbb{P}^n$-bundles.

Harder-Katzarkov define a version of a perverse sheaf of categories motivated by the notion of perverse Schobers due to Kapranov-Schechtman \cite{kapranovschechtman}.
These notions are quiver theoretic in nature, and following this line the CCC for principal toric fibrations can be studied using quivers.
This is well expected in that calculation of microlocal sheaf categories can be done in a combinatorial way \cite{arboreal},
and such categories have nice quiver models when the skeleta resemble combinatorial nature,
which is indeed the case for FLTZ skeleta.

For example, consider $\ush_{\Lambda_{\mathbb{P}^1}}$.
This sheaf of categories is supported on $\Lambda_{\mathbb{P}^1}$ and, in fact, constructible on $\Lambda_{\mathbb{P}^1}$. That is, considering $\Lambda_{\mathbb{P}^1}$ as a stratified space (whose strata are a point, an open interval, and two half rays), $\ush_{\Lambda_{\mathbb{P}^1}}$ is locally constant (in fact constant) on these contractible strata. Hence to describe $\ush_{\Lambda_{\mathbb{P}^1}}$, we need the following information: four categories that correspond to stalks of $\ush_{\Lambda_{\mathbb{P}^1}}$ at each stratum and co-restriction morphisms
for each inclusion of strata.

The case of $\ush_{\Lambda_{\mathbb{P}^1}} \otimes_{k} \Phi(\mathcal{E})$ is similar,
with the exception that there is a monodromy induced from $\Phi(\mathcal{E})$.
Harder-Katzarkov incorporates this monodromy by twisting morphisms between categories, leading to Definition \ref{def:HK}.

Taking the global section of $\ush_{\Lambda_{\mathbb{P}^1}}$, we obtain the category of representations of the followin gquiver, which is exactly Beilinson's quiver for $\mathrm{D}^b(\mathbb{P}^2)$:
\[
b \leftarrow a \rightarrow b
\]

For the toric $\mathbb{P}^1$-bundle $\mathbb{P}(\mathcal{O}\oplus \mathcal{L})$, $\Phi(\mathbb{P}(\mathcal{O}\oplus \mathcal{L}))$ is a local system whose monodromy around $S^1$ is $-\otimes \mathcal{L}$. Hence we expect that the global section of $\ush_{\Lambda_{\mathbb{P}^1}} \otimes_{k} \Phi(\mathbb{P}(\mathcal{O}\oplus \mathcal{L}))$ is the category of representations of the following quiver,
\[
b \leftarrow a \rightarrow b\otimes \mathcal{L},
\]
and this coincides with the one considered by Harder-Katzarkov in the context of perverse sheaves of categories.
The theorem of Harder-Katzarkov (Theorem \ref{thm:HK}) asserts that this category is actually equivalent to $\mathrm{D}^b(\mathbb{P}(\Oc\oplus \Lc))$.

We can interpret this as follows: by the monodromy $-\otimes \mathcal{L}$, the left $b$ varies to $b\otimes\mathcal{L}$ as we go around $S^1$.

\begin{figure}\label{fig:constrP2}
	\begin{tikzpicture}
\newcommand*{\edgelen}{4}; \newcommand*{\vertrad}{.12}; \newcommand*{\crad}{.05}
	\newcommand*{\gspace}{.9}; \newcommand*{\angdelta}{22.5};

	\node (a) [matrix] at (-3,0) {
		\coordinate (tl) at (0,\edgelen); \coordinate (tr) at (\edgelen,\edgelen);
		\coordinate (bl) at (0,0); \coordinate (br) at (\edgelen,0); 
		\coordinate (tm) at (\edgelen/3,\edgelen); \coordinate (bm) at (\edgelen/3,0);
		\coordinate (mm) at (\edgelen/3,\edgelen/3);
		\coordinate (lm) at (0, \edgelen/3); \coordinate (rm) at (\edgelen, \edgelen/3);
		\coordinate (lmin) at (\edgelen/3, \edgelen*2/3); \coordinate (rmin) at (\edgelen*2/3, \edgelen/3);
		\foreach \a/\b in {tl/bl, tr/br} {\draw[asdstyle] (\a) to (\b);}
		\foreach \a/\b in {tr/tl, br/bl} {\draw[asdstyle] (\a) to (\b);}
		\foreach \a/\b in {lm/rm} {\draw[asdstyle,righthairs] (\a) to (\b);}
		\foreach \a/\b in {bm/tm, tl/br} {\draw[asdstyle,lefthairs] (\a) to (\b);}
		\foreach \c in {mm} \foreach \ang in {8,...,12}
			{\draw[blue,thick] ($(\c)+(\ang*\angdelta:\vertrad+.02)$) to (\c);}
		\foreach \c in {lmin} \foreach \ang in {2,...,8}
			{\draw[blue,thick] ($(\c)+(\ang*\angdelta:\vertrad+.02)$) to (\c);}
		\foreach \c in {rmin} \foreach \ang in {0,1,2,12,13,14,15,16}
			{\draw[blue,thick] ($(\c)+(\ang*\angdelta:\vertrad+.02)$) to (\c);}
		\node (tltri) at (\edgelen*2/9, \edgelen*7/8) {$a$};
		\node (square) at (\edgelen*1.5/9, \edgelen*1.5/9) {$a$};
		\node (brtri) at (\edgelen*8.2/9, \edgelen*1.5/9) {$a$};
		\node (hexa) at (\edgelen*2/3, \edgelen*2/3) {$b$};
		\node (btrap) at (\edgelen*1/2, \edgelen*1.5/9) {$b$};
		\node (ttrap) at (\edgelen*1.5/9, \edgelen*1/2) {$b$};
		\node (center) at (\edgelen*3/7, \edgelen*3/7) {$c$};
		\draw[genmapstyle] (tltri) to (ttrap);
		\draw[genmapstyle] (tltri) to (hexa);
		\draw[genmapstyle] (square) to (ttrap);
		\draw[genmapstyle] (square) to (btrap);
		\draw[genmapstyle] (brtri) to (btrap);
		\draw[genmapstyle] (brtri) to (hexa);
		\draw[genmapstyle] (ttrap) to (center);
		\draw[genmapstyle] (btrap) to (center);
		\draw[genmapstyle] (hexa) to (center);
	\\};

	\node (b) [matrix] at (3, 0) {
		\coordinate (tl) at (0,\edgelen); \coordinate (tr) at (\edgelen,\edgelen);
		\coordinate (bl) at (0,0); \coordinate (br) at (\edgelen,0); 
		\coordinate (tm) at (\edgelen/3,\edgelen); \coordinate (bm) at (\edgelen/3,0);
		\coordinate (mm) at (\edgelen/3,\edgelen/3);
		\coordinate (lm) at (0, \edgelen/3); \coordinate (rm) at (\edgelen, \edgelen/3);
		\coordinate (lmin) at (\edgelen/3, \edgelen*2/3); \coordinate (rmin) at (\edgelen*2/3, \edgelen/3);
		\foreach \a/\b in {tl/bl, tr/br} {\draw[asdstyle] (\a) to (\b);}
		\foreach \a/\b in {tr/tl, br/bl} {\draw[asdstyle] (\a) to (\b);}
		\foreach \a/\b in {lm/rm} {\draw[asdstyle,righthairs] (\a) to (\b);}
		\foreach \a/\b in {bm/tm, tl/br} {\draw[asdstyle,lefthairs] (\a) to (\b);}
		\foreach \c in {mm} \foreach \ang in {8,...,12}
			{\draw[blue,thick] ($(\c)+(\ang*\angdelta:\vertrad+.02)$) to (\c);}
		\foreach \c in {lmin} \foreach \ang in {2,...,8}
			{\draw[blue,thick] ($(\c)+(\ang*\angdelta:\vertrad+.02)$) to (\c);}
		\foreach \c in {rmin} \foreach \ang in {0,1,2,12,13,14,15,16}
			{\draw[blue,thick] ($(\c)+(\ang*\angdelta:\vertrad+.02)$) to (\c);}
		\node (tltri) at (\edgelen*2/9, \edgelen*7/8) {$a\mathcal{M}$};
		\node (square) at (\edgelen*1.5/9, \edgelen*1.5/9) {$a$};
		\node (brtri) at (\edgelen*8.2/9, \edgelen*1.5/9) {$a\mathcal{L}$};
		\node (hexa) at (\edgelen*2/3, \edgelen*2/3) {$b\mathcal{M}$};
		\node (btrap) at (\edgelen*1/2, \edgelen*1.5/9) {$b$};
		\node (ttrap) at (\edgelen*1.5/9, \edgelen*1/2) {$b\mathcal{L}^{-1}\mathcal{M}$};
		\node (center) at (\edgelen*3/7, \edgelen*3/7) {$c$};
		\draw[genmapstyle] (tltri) to (ttrap);
		\draw[genmapstyle] (tltri) to (hexa);
		\draw[genmapstyle] (square) to (ttrap);
		\draw[genmapstyle] (square) to (btrap);
		\draw[genmapstyle] (brtri) to (btrap);
		\draw[genmapstyle] (brtri) to (hexa);
		\draw[genmapstyle] (ttrap) to (center);
		\draw[genmapstyle] (btrap) to (center);
		\draw[genmapstyle] (hexa) to (center);
	\\};

	\end{tikzpicture}
	\caption{$\ush_{\Lambda_{\mathbb{P}^2}}$ vs $\ush_{\Lambda_{\mathbb{P}^2}} \otimes \Phi(\mathbb{P}(\mathcal{O}\oplus \mathcal{L} \oplus \mathcal{M}))$}\label{comparison}
\end{figure}
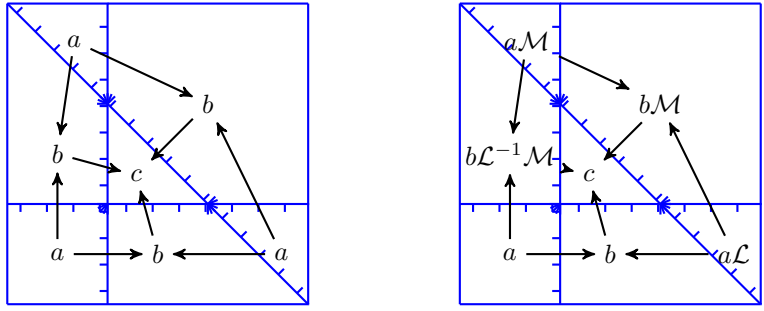

Although the full theory of perverse schobers on higher dimensional spaces is not developed yet,
similar to the above case of $\mathbb{P}(\mathcal{O}\oplus \mathcal{L})$,
for $\mathbb{P}(\mathcal{O}\oplus \mathcal{L} \oplus \mathcal{M})$,
we may expect that the global section of $\ush_{\Lambda_{\mathbb{P}^2}}$ and $\ush_{\Lambda_{\mathbb{P}^2}} \otimes_{k} \Phi(\mathbb{P}(\mathcal{O}\oplus \mathcal{L} \oplus \mathcal{M}))$ can be expressed in terms of the category of quiver representations in the Figure \ref{comparison}.
In Figure \ref{comparison}, we describe $T^2$ by the blue square with each side identified, and the whole blue picture denotes the (perturbation of) $\Lambda_{\mathbb{P}^2}$. Now $\Phi(\mathbb{P}(\mathcal{O}\oplus \mathcal{L} \oplus \mathcal{M}))$ is a local system on $T^*T^2$ whose monodromies are $-\otimes \mathcal{L}$ and $-\otimes \mathcal{M}$. Hence as we go around along the horizontal $S^1$, $a$ varies to $a\mathcal{L}$. Along the vertical $S^1$, $a$ varies to $a\mathcal{M}$. Similarly, we have $b\mathcal{M}$ and $b\mathcal{L}^{-1}\mathcal{M}$.

As an analog of Theorem \ref{thm:HK}, one expects that the category of such quiver representations is equivalent to $\mathrm{D}^b(\mathbb{P}(\Oc\oplus \Lc \oplus \Mc))$. In fact, we can prove this (Theorem \ref{thm:P2}) and generalize to the case of any toric $\mathbb{P}^n$-bundles(Theorem \ref{thm:Pn}) independent of our main theorem, following the idea of Harder-Katzarkov's proof of Theorem \ref{thm:HK}.
In the remaining of this chapter, we recall Harder-Katzarkov's result in detail and generalize their proof to
the case of $\mathbb{P}^n$-bundles.

\subsection{Perverse schobers and perverse sheaf of categories}

In \cite{harder2019perverse}, the authors define a perverse sheaf of categories on a surface as a generalization of perverse schobers, as categorifications of perverse sheaves \cite{kapranovschechtman}. 

Recall that the category $\mathrm{Perv}(D,0)$ of perverse sheaves on a unit disk $D$ in $\mathbb{C}$ stratified by the origin and its complement admits a quiver description.

\begin{theorem}\cite{GGM}
$\mathrm{Perv}(D, 0)$ is equivalent to the category of quadruples $(V,W,\phi, \psi)$ where $V,W$ are vector spaces and $\phi:V\to W$, $\psi:W \to V$ are linear maps such that $\id_V-\psi \circ \phi$ and $\id_W-\phi \circ \psi$ are invertible. 
\end{theorem}

To get an idea why this is true, fix a point $b$ on the boundary of $D$, and choose an arc $K$ connecting the origin and $b$. Note that $K$ is a skeleton of $D$.
Then to $\mathcal{F} \in \mathrm{Perv}(D, 0)$, we associate $V(\mathcal{F})=\mathbb{H}^1_{K}(\mathcal{F})_0$ and $W(\mathcal{F})=\mathbb{H}^1_{K}(\mathcal{F})_b$. Here, $\mathbb{H}^1_K$ denotes the first hypercohomology with support in $K$, and $\phi$ is the corestriction map for the constructible sheaf $\mathbb{H}^1_{K}(\mathcal{F})$, while $\psi$ is the composition of the counterclockwise continuation of $\mathcal{F}_b$ and the differential 
\[
H^{0}(\Delta - K, \mathcal{F}) \to \mathbb{H}^1_{K}(\Delta, \mathcal{F}).
\]
It turns out this procedure gives an equivalence of categories.

Motivated by this theorem, \cite{kapranovschechtman}
suggests a natural categorification of a perverse sheaf on $D$ be the notion of a spherical functor:

\begin{definition}[\cite{annologvinenko}]
Let $S:\mathcal{C} \to \mathcal{D}$ be an exact functor between dg triangulated categories $\mathcal{C}$ and $\mathcal{D}$. Assume that $S$ admits left and right adjoints $L$ and $R$. Then we have four natural transformations 
\[SR \to \id_{\mathcal{D}}, LS \to \id_{\mathcal{C}}\]
\[\id_{\mathcal{C}} \to RS, \id_{\mathcal{D}} \to SL.\]

We denote the cones of these natural transformations by 
\[T=\mathrm{Cone}(SR \to \id_{\mathcal{D}}), C'=\mathrm{Cone}(LS \to \id_{\mathcal{C}})\]
\[C=\mathrm{Cone}(\id_{\mathcal{C}} \to RS)[-1], T'=\mathrm{Cone}(\id_{\mathcal{D}} \to SL)[-1].\]
$S$ is called a \emph{spherical functor} if the following hold\footnote{In fact, any two of them imply the other two, see \cite{annologvinenko}.
}:
\begin{enumerate}
    \item $T$ is an equivalence.
    \item The composition $R\to RSL \to CL[1]$ is an isomorphism.
    \item $C$ is an equivalence.
    \item The composition $LT[-1]\to LSR \to R$ is an isomorphism.
\end{enumerate}
\end{definition}
Similar results continue to hold in the case of multiple marked points.
Let $\Sigma=\{p_1, \cdots, p_n\}$ be a set of marked points in $D$.
The category $\mathrm{Perv}(D, \Sigma)$ of perverse sheaves on $D$ stratified by $\Sigma$ and its complement admits the following description.

\begin{theorem}[\cite{GMV}]
$\mathrm{Perv}(D, \Sigma)$ is equivalent to the category of diagrams formed by vector spaces $V, W_1, \cdots, W_n$ and linear maps $\phi_i:V\to W_i$, $\psi_i:W_i \to V$ such that $\id_V-\psi_i \circ \phi_i$ and $\id_{W_i}-\phi_i\circ\psi_i$ are invertible. 
\end{theorem}

This description of $\mathrm{Perv}(D, \Sigma)$ motivates
the following definition.

\begin{definition}
A $K$-coordinatized perverse schober on $(D, \Sigma)$ is a system of categories $\mathcal{C}, \mathcal{D}_1, \cdots, \mathcal{D}_n$ together with spherical functors $S_i : \mathcal{C} \to \mathcal{D}_i$.
\end{definition}

\subsection{CCC for toric $\mathbb{P}^1$-bundles due to Harder-Katzarkov}

In this section, we review the theorem of Harder-Katzarkov on CCC for toric $\mathbb{P}^1$-bundles. A reference for this section is \cite{harder2019perverse}. 

Motivated by the definition of perverse schobers, \cite{harder2019perverse} define a perverse sheaf of categories on a surface $S$ as follows.

Let $S$ be a real $2$-dimensional compact, connected, oriented surface with $n$ boundary components and $k$ marked points $\Sigma$.
\begin{definition}
A graph $K$ inside the surface $S$ is a collection of
\begin{itemize}
\item a finite subset $Vert(K)$ of the interior of $S$,
\item a finite set $Edge(K)$ of embedded closed intervals in $S$,
with ends being either a vertex or on the boundary of $S$,
while at least one of the ends is a vertex,
\end{itemize}
subject to the further assumptions that the interior of an edge lies in the interior of $S$ and any two edges only intersect at vertices of $K$.
\end{definition}

\begin{construction}\label{An-perverse-schober}
Let $A_n(\mathcal{C})$ denote the category of representations of the $A_n$-quiver over $\mathcal{C}$. It comes with $n+1$ natural functors which we denote by $f_1, \cdots, f_{n+1}$ from $A_n(\mathcal{C}) \to \mathcal{C}$. For example, if $n=2$, an object of $A_2(\mathcal{C})$ is roughly speaking $x \xrightarrow{a} y$ for $x,y \in ob(\mathcal{C})$ and $a\in \mathrm{Hom}_{\mathcal{C}}(x, y)$. Then $f_i : A_2(\mathcal{C}) \to \mathcal{C}$ are given by $f_1(x \xrightarrow{a} y) = x, f_2(x \xrightarrow{a} y) = y, f_3(x \xrightarrow{a} y) = \mathrm{Cone}(a)$. The detailed construction of $A_n(\mathcal{C})$ and the $f_i$ can be found in \cite{harder2019perverse}.
\end{construction}
\begin{definition}
A $K$-coordinatized perverse sheaf of categories on $(S,\Sigma)$ is a collection of data
consisting of a choice of a skeleton $K$ of $S$, a pretriangulated dg category $\mathcal{C}$, a collection of pretriangulated dg categories $\mathcal{A}_v$ for all $v\in Vert(K)$ and functors $F_{v,e}:\mathcal{A}_v \to \mathcal{C}$ for all edgs $e\in Edge(K)$,  together with
a monodromy representation
$\pi_1(S^{\circ}, v_1) \rightarrow \operatorname{Aut}(H^0 \mathcal{C})$,
satisfying the following conditions:
\begin{enumerate}
\item $K$ is a spanning graph of $S$ i.e. $K$ is a graph inside $S$ that is homotopic to $S$ itself. We further require that any marked point is a univalent vertex of $K$ and that there are no bivalent vertices.

\item A category $\mathcal{A}_{v}$ for each vertex $v$ in $Vert(K)$ and a fixed category $\mathcal{C}$ for each edge $e$ in $Edge(K)$. If $v$ is not in $\Sigma$, then we assign $\mathcal{A}_v=A_n(\mathcal{C})$, where $n+1$ is the valency of $v$.

\item If an edge $e$ is incident to a vertex $v$, we assign a functor $F_{v,e}:\mathcal{A}_v \to \mathcal{C}$. If $v$ is $n+1$-valent, then for some counterclockwise ordering of edges incident to $v$, denoted $e_1, \cdots e_{n+1}$, we have $F_{v, e_i}=\phi_{v, e_i} \circ f_i$ for some autoequivalence $\phi_{v, e_i}$ of $\mathcal{C}$. If $v\in\Sigma$, then we assume that $F_{v,e}$ is spherical.
\end{enumerate}
\label{def:HK}\end{definition}

A $K$-coordinatized perverse sheaf of categories $\mathcal{F}$ induces a diagram of categories, given by assigning the category $\mathcal{A}_v$ to each vertex $v$ and functors $F_{v,e}$ to each edge $e$.
The homotopy limit of this diagram, denoted by $\Gamma(K, \mathcal{F})$ is called the global section of $\mathcal{F}$.

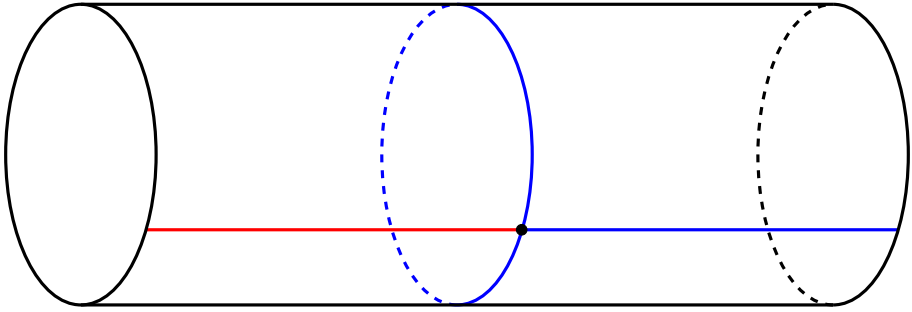
\begin{figure}
    \begin{tikzpicture}

\draw[fill = black] (0.86,-1) circle (1pt);

\draw[blue,very thick] (0.86,-1) to (5.86,-1);
\draw[red,very thick] (0.86,-1) to (-4.14,-1);

	\draw[very thick] (5,-2) arc (270:450:1cm and 2cm);
    \draw[dashed,very thick] (5,-2) arc (270:90:1cm and 2cm);
	\draw[very thick] (-5,0) ellipse (1cm and 2cm);

\draw[blue,very thick] (0,-2) arc (270:450:1cm and 2cm);
\draw[dashed,blue, very thick] (0,-2) arc (270:90:1cm and 2cm);
\draw[very thick] (-5,2) to (5,2);
\draw[very thick] (-5,-2) to (5,-2);

\draw[fill = black] (0.86,-1) circle (2pt);

\end{tikzpicture}\caption{\label{P1skel}$\Lambda_{\Sigma_{\mathbb{P}^1}}$}

\end{figure}

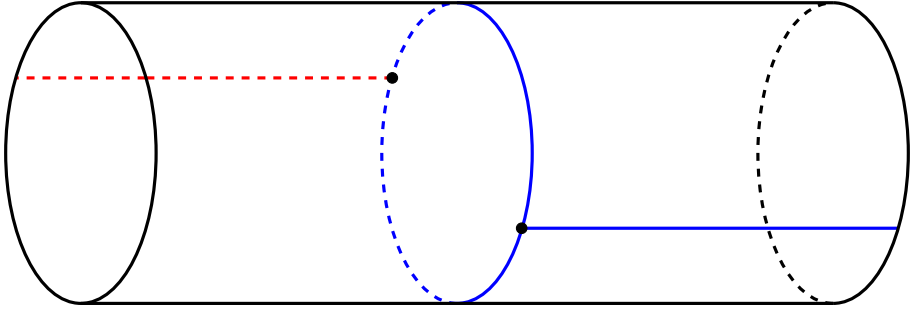
\begin{figure}
    \begin{tikzpicture}

\draw[fill = black] (0.86,-1) circle (1pt);
\draw[fill = black] (-0.86,1) circle (1pt);

\draw[blue,very thick] (0.86,-1) to (5.86,-1);
\draw[red,dashed,very thick] (-0.86,1) to (-4.14,1);
\draw[red, very thick,dashed] (-4.14,1) to (-5.86,1);

	\draw[very thick] (5,-2) arc (270:450:1cm and 2cm);
    \draw[dashed,very thick] (5,-2) arc (270:90:1cm and 2cm);
	\draw[very thick] (-5,0) ellipse (1cm and 2cm);

\draw[blue,very thick] (0,-2) arc (270:450:1cm and 2cm);
\draw[dashed,blue, very thick] (0,-2) arc (270:90:1cm and 2cm);
\draw[very thick] (-5,2) to (5,2);
\draw[very thick] (-5,-2) to (5,-2);

\draw[fill = black] (0.86,-1) circle (2pt);
\draw[fill = black] (-0.86,1) circle (2pt);
\end{tikzpicture}\caption{\label{P1skelprime}$\widetilde{\Lambda}_{\Sigma_{\mathbb{P}^1}}$}

\end{figure}

\begin{example}
Let $\Lambda_{\Sigma_{\mathbb{P}^1}}$ be the conical Lagrangian skeleton of $T^*T^1$ associated to $\mathbb{P}^1$ under CCC  (Figure \ref{P1skel}) and $\widetilde{\Lambda}_{\Sigma_{\mathbb{P}^1}}$ be its perturbation depicted in Figure \ref{P1skelprime}.\footnote{Figures credit to \cite{harder2019perverse}.} Let $\mathcal{L}$ be a line bundle on a variety $X$.
We define a perverse sheaf of categories $\mathcal{F}_\mathcal{L}$ on $\widetilde{\Lambda}_{\Sigma_{\mathbb{P}^1}}$ as follows. Let $v_1$ be the upper vertex and $v_2$ be the lower vertex in Figure \ref{P1skelprime}. Also, we denote the upper and lower half circles connecting $v_1$ and $v_2$ by $e_1$ and $e_2$ respectively. Let $e_3$ be the red straight line, and $e_4$ be the blue straight line. We set $\mathcal{A}_{v_1}=\mathcal{A}_{v_2}=A_2(\mathcal{C})$,
and define the functors from vertices to edges by $F_{v_1, e_1}=F_{v_2, e_1} =f_1, F_{v_1, e_2}=f_2, F_{v_2, e_2} = (-\otimes \mathcal{L}) \circ f_2, F_{v_1, e_3}=F_{v_2, e_4}=f_3$.\footnote{Recall Construction \ref{An-perverse-schober}.}

\end{example}

It's evident that $\mathbb{P}(\mathcal{O} \oplus \mathcal{L})$ is a
principal toric fibration with fiber $\mathbb{P}^1$ over $X$.
Let $\mathcal{C}=\mathrm{D}^b(X)$ in the above example. Harder-Katzarkov prove the following theorem.

\begin{theorem}{(CCC for a toric $\mathbb{P}^1$- bundle)}\cite[Theorem 4.5]{harder2019perverse}
The dg categories $\mathrm{D}^b(\mathbb{P}(\mathcal{O} \oplus \mathcal{L}))$ and $\Gamma(\widetilde{\Lambda}_{\Sigma_{\mathbb{P}^1}}, \mathcal{F}_{\mathcal{L}})$ are quasi-equivalent.
\label{thm:HK}
\end{theorem}

\subsection{Quiver theoretic approach to CCC for toric $\mathbb{P}^n$-bundles}

As we can see, Harder-Katzarkov's approach is quiver-theoretic in nature.
The perverse sheaf of categories
$\Fc_\Lc$
they defined can be seen as a concrete model for 
$\ush(\Lambda_{\mathbb{P}^1}) \otimes \Phi(\mathbb{P}(\mathcal{O} \oplus \mathcal{L}))$. $\Gamma(\widetilde{\Lambda}_{\Sigma_{\mathbb{P}^1}}, \mathcal{F}_{\mathcal{L}})$ is a subcategory of the category of representations of the quiver $\cdot \leftarrow \cdot \rightarrow \cdot$ with coefficients in $\mathrm{D}^b(X)$ consisting of objects of the form $a \leftarrow b \rightarrow a \otimes \mathcal{L}$ for $a,b \in ob(\mathrm{D}^b(X))$.

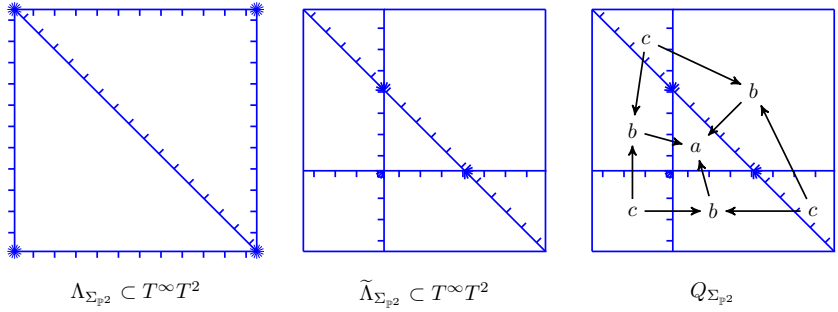
\begin{figure}
    \resizebox{0.9\textwidth}{!}{
	\begin{tikzpicture}[scale=0.9]
	\newcommand*{\edgelen}{4}; \newcommand*{\vertrad}{.12}; \newcommand*{\crad}{.05}
	\newcommand*{\gspace}{.9}; \newcommand*{\angdelta}{22.5};
	\node (a) [matrix] at (-5.3,0) {
		\coordinate (tl) at (0,\edgelen); \coordinate (tr) at (\edgelen,\edgelen);
		\coordinate (bl) at (0,0); \coordinate (br) at (\edgelen,0);
		\foreach \a/\b in {tl/bl, tr/br, br/tl} {\draw[asdstyle,righthairs] (\a) to (\b);}
		\foreach \a/\b in {tr/tl, br/bl} {\draw[asdstyle,lefthairs] (\a) to (\b);}
		\foreach \c in {tl, tr, bl, br} \foreach \ang in {0,...,7}
			{\draw[blue] ($(\c)+(\ang*\angdelta:\vertrad)$) to ($(\c)+(\ang*\angdelta+180:\vertrad)$);}
	\\};
	\node (b) [matrix] at (0,0) {
		\coordinate (tl) at (0,\edgelen); \coordinate (tr) at (\edgelen,\edgelen);
		\coordinate (bl) at (0,0); \coordinate (br) at (\edgelen,0); 
		\coordinate (tm) at (\edgelen/3,\edgelen); \coordinate (bm) at (\edgelen/3,0);
		\coordinate (mm) at (\edgelen/3,\edgelen/3);
		\coordinate (lm) at (0, \edgelen/3); \coordinate (rm) at (\edgelen, \edgelen/3);
		\coordinate (lmin) at (\edgelen/3, \edgelen*2/3); \coordinate (rmin) at (\edgelen*2/3, \edgelen/3);
		\foreach \a/\b in {tl/bl, tr/br} {\draw[asdstyle] (\a) to (\b);}
		\foreach \a/\b in {tr/tl, br/bl} {\draw[asdstyle] (\a) to (\b);}
		\foreach \a/\b in {lm/rm} {\draw[asdstyle,righthairs] (\a) to (\b);}
		\foreach \a/\b in {bm/tm, tl/br} {\draw[asdstyle,lefthairs] (\a) to (\b);}
		\foreach \c in {mm} \foreach \ang in {8,...,12}
			{\draw[blue,thick] ($(\c)+(\ang*\angdelta:\vertrad+.02)$) to (\c);}
		\foreach \c in {lmin} \foreach \ang in {2,...,8}
			{\draw[blue,thick] ($(\c)+(\ang*\angdelta:\vertrad+.02)$) to (\c);}
		\foreach \c in {rmin} \foreach \ang in {0,1,2,12,13,14,15,16}
			{\draw[blue,thick] ($(\c)+(\ang*\angdelta:\vertrad+.02)$) to (\c);}
	\\};
	
	\node (a) [matrix] at (5.3,0) {
		\coordinate (tl) at (0,\edgelen); \coordinate (tr) at (\edgelen,\edgelen);
		\coordinate (bl) at (0,0); \coordinate (br) at (\edgelen,0); 
		\coordinate (tm) at (\edgelen/3,\edgelen); \coordinate (bm) at (\edgelen/3,0);
		\coordinate (mm) at (\edgelen/3,\edgelen/3);
		\coordinate (lm) at (0, \edgelen/3); \coordinate (rm) at (\edgelen, \edgelen/3);
		\coordinate (lmin) at (\edgelen/3, \edgelen*2/3); \coordinate (rmin) at (\edgelen*2/3, \edgelen/3);
		\foreach \a/\b in {tl/bl, tr/br} {\draw[asdstyle] (\a) to (\b);}
		\foreach \a/\b in {tr/tl, br/bl} {\draw[asdstyle] (\a) to (\b);}
		\foreach \a/\b in {lm/rm} {\draw[asdstyle,righthairs] (\a) to (\b);}
		\foreach \a/\b in {bm/tm, tl/br} {\draw[asdstyle,lefthairs] (\a) to (\b);}
		\foreach \c in {mm} \foreach \ang in {8,...,12}
			{\draw[blue,thick] ($(\c)+(\ang*\angdelta:\vertrad+.02)$) to (\c);}
		\foreach \c in {lmin} \foreach \ang in {2,...,8}
			{\draw[blue,thick] ($(\c)+(\ang*\angdelta:\vertrad+.02)$) to (\c);}
		\foreach \c in {rmin} \foreach \ang in {0,1,2,12,13,14,15,16}
			{\draw[blue,thick] ($(\c)+(\ang*\angdelta:\vertrad+.02)$) to (\c);}
		\node (tltri) at (\edgelen*2/9, \edgelen*7/8) {$c$};
		\node (square) at (\edgelen*1.5/9, \edgelen*1.5/9) {$c$};
		\node (brtri) at (\edgelen*8.2/9, \edgelen*1.5/9) {$c$};
		\node (hexa) at (\edgelen*2/3, \edgelen*2/3) {$b$};
		\node (btrap) at (\edgelen*1/2, \edgelen*1.5/9) {$b$};
		\node (ttrap) at (\edgelen*1.5/9, \edgelen*1/2) {$b$};
		\node (center) at (\edgelen*3/7, \edgelen*3/7) {$a$};
		\draw[genmapstyle] (tltri) to (ttrap);
		\draw[genmapstyle] (tltri) to (hexa);
		\draw[genmapstyle] (square) to (ttrap);
		\draw[genmapstyle] (square) to (btrap);
		\draw[genmapstyle] (brtri) to (btrap);
		\draw[genmapstyle] (brtri) to (hexa);
		\draw[genmapstyle] (ttrap) to (center);
		\draw[genmapstyle] (btrap) to (center);
		\draw[genmapstyle] (hexa) to (center);
	\\};

	\node (atxt) at (-5.3, -3) {$\Lambda_{\Sigma_{\mathbb{P}^2}} \subset T^\infty T^2$};
	\node (btxt) at (0, -3) {$\widetilde{\Lambda}_{{\Sigma}_{\mathbb{P}^2}} \subset T^\infty T^2$};
	
	\node (ctxt) at (5.3, -3)
	{$Q_{\Sigma_{\mathbb{P}^2}}$};
	\end{tikzpicture}
}
	\caption{The FLTZ skeleton of $\mathbb{P}^2$}\label{fig:fltzP2}
\end{figure}

\begin{figure}
\resizebox{0.7\textwidth}{!}{
\begin{tikzpicture}[scale=0.7]
\newcommand*{\edgelen}{12}; \newcommand*{\vertrad}{.12}; \newcommand*{\crad}{.05}
	\newcommand*{\gspace}{.9}; \newcommand*{\angdelta}{22.5};

	\node (c) [matrix] at (5.3,0) {
		\coordinate (tl) at (0,\edgelen); \coordinate (tr) at (\edgelen,\edgelen);
		\coordinate (bl) at (0,0); \coordinate (br) at (\edgelen,0); 
		\coordinate (tm) at (\edgelen/3,\edgelen); \coordinate (bm) at (\edgelen/3,0);
		\coordinate (mm) at (\edgelen/3,\edgelen/3);
		\coordinate (lm) at (0, \edgelen/3); \coordinate (rm) at (\edgelen, \edgelen/3);
		\coordinate (lmin) at (\edgelen/3, \edgelen*2/3); \coordinate (rmin) at (\edgelen*2/3, \edgelen/3);
		\foreach \a/\b in {tl/bl, tr/br} {\draw[asdstyle] (\a) to (\b);}
		\foreach \a/\b in {tr/tl, br/bl} {\draw[asdstyle] (\a) to (\b);}
		\foreach \a/\b in {lm/rm} {\draw[asdstyle,righthairs] (\a) to (\b);}
		\foreach \a/\b in {bm/tm, tl/br} {\draw[asdstyle,lefthairs] (\a) to (\b);}
		\foreach \c in {mm} \foreach \ang in {8,...,12}
			{\draw[blue,thick] ($(\c)+(\ang*\angdelta:\vertrad+.02)$) to (\c);}
		\foreach \c in {lmin} \foreach \ang in {2,...,8}
			{\draw[blue,thick] ($(\c)+(\ang*\angdelta:\vertrad+.02)$) to (\c);}
		\foreach \c in {rmin} \foreach \ang in {0,1,2,12,13,14,15,16}
			{\draw[blue,thick] ($(\c)+(\ang*\angdelta:\vertrad+.02)$) to (\c);}
		\node (tltri) at (\edgelen*2/9, \edgelen*7/8) {$x\mathcal{M}$};
		\node (square) at (\edgelen*1.5/9, \edgelen*1.5/9) {$x$};
		\node (brtri) at (\edgelen*8.2/9, \edgelen*1.5/9) {$x\mathcal{L}^{-1}\mathcal{M}$};
		\node (hexa) at (\edgelen*2/3, \edgelen*2/3) {$y\mathcal{M}$};
		\node (btrap) at (\edgelen*1/2, \edgelen*1.5/9) {$y$};
		\node (ttrap) at (\edgelen*1.5/9, \edgelen*1/2) {$y\mathcal{L}$};
		\node (center) at (\edgelen*3/7, \edgelen*3/7) {$z$};
		\draw[genmapstyle] (tltri) -- node [left] {$h\mathcal{L}$} (ttrap);
		\draw[genmapstyle] (tltri) -- node [above] {$f\mathcal{M}$} (hexa);
		\draw[genmapstyle] (square) --node [left] {$g$} (ttrap);
		\draw[genmapstyle] (square) -- node [above] {$f$} (btrap);
		\draw[genmapstyle] (brtri) -- node [above] {$h$} (btrap);
		\draw[genmapstyle] (brtri) -- node [right] {$g\mathcal{L}^{-1}\mathcal{M}$} (hexa);
		\draw[genmapstyle] (ttrap) -- node [above] {$k_1$} (center);
		\draw[genmapstyle] (btrap) -- node [left] {$k_3$} (center);
		\draw[genmapstyle] (hexa) -- node [above] {$k_2$} (center);
	\\};

	\end{tikzpicture}
}
    \caption{An object of $Rep(Q_{\Sigma_{\mathbb{P}^n}}, \mathrm{D}^b(X), \{\mathcal{L}, \mathcal{M}\})$}
    \label{fig:typical}
\end{figure}
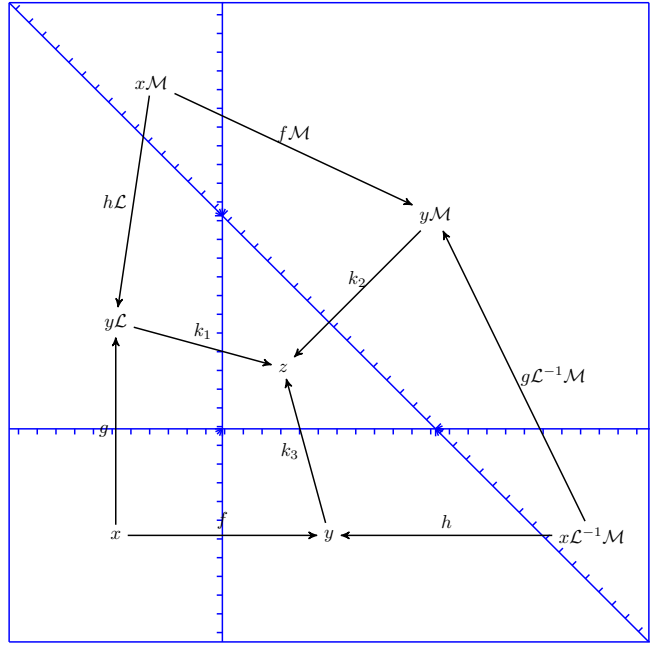

Their approach indeed is applicable to the case of other toric varieties. That is, we may express $\Gamma(\Lambda_\Sigma ,\ush(\Lambda_\Sigma)\otimes \Phi\mathcal(X_{\Sigma}))$ in terms of quiver representations.
Let $\Lambda_{\Sigma_{\mathbb{P}^2}}$ be the FLTZ skeleton for $\mathbb{P}^2$
which is shown on the left of Figure \ref{fig:fltzP2}.
To describe objects in this category, for the purpose of easy visualization, we first perturb $\Lambda_{\Sigma_{\mathbb{P}^2}}$ into $\widetilde{\Lambda}_{\Sigma_{\mathbb{P}^2}}$ as in Figure \ref{fig:fltzP2}. Now, an object of this category of $\mathbb{P}^2$ can be expressed as in Figure \ref{fig:typical}.
Similarly an object in $\Gamma(T^*T^2, \ush_{\Lambda_{\Sigma_{\mathbb{P}^2}}})$ can be described by the third figure in Figure \ref{fig:fltzP2} where $a,b,c \in D^b_{dg}(pt)=Vect$.
Note that this is exactly a representation of Beilinson's quiver for $\mathbb{P}^2$, albeit with
coefficients in the category $\mathrm{D}^b(X)$
together with an extra twist.

We claim that a corresponding category for $\mathbb{P}(\mathcal{O} \oplus \mathcal{L} \oplus \mathcal{M})$, which we denote by $Rep(Q_{\Sigma_{\mathbb{P}^n}}, \mathrm{D}^b(X), \{\mathcal{L}, \mathcal{M}\})$ is the following. It is a full subcategory of the category of representations of the same quiver for $\mathbb{P}^2$, which we denote by $Q_{\Sigma_{\mathbb{P}^2}}$, with coefficients in $\mathrm{D}^b(X)$ consisting of objects of the form described in Figure \ref{fig:typical}.

Note that if $X$ is a point, then $\mathbb{P}(\mathcal{O} \oplus \mathcal{L} \oplus \mathcal{M})$ is $\mathbb{P}^2$. In this case, $$Rep(Q_{\Sigma_{\mathbb{P}^n}}, \mathrm{D}^b(X), \{\mathcal{L}, \mathcal{M}\})$$ is just the category of representation of Beilinson's quiver for $\mathbb{P}^2$.

For toric $\mathbb{P}^2$-bundle, our main theorem can be re-stated as follows.

\begin{theorem}
$\mathrm{D}^b(\mathbb{P}(\mathcal{O}\oplus \mathcal{L} \oplus \mathcal{M}))$ is quasi-equivalent to
$Rep(Q_{\Sigma_{\mathbb{P}^n}}, \mathrm{D}^b(X), \{\mathcal{L}, \mathcal{M}\})$.
\label{thm:P2}
\end{theorem}

Let us give another proof of this theorem, in a similar way to Harder-Katzarkov's proof for Theorem \ref{thm:HK}.

\begin{proof}
We first prove that $\mathrm{D}^b(\mathbb{P}(\mathcal{O}\oplus \mathcal{L} \oplus \mathcal{M}))$ embeds into $Rep(Q_{\Sigma_{\mathbb{P}^n}}, \mathrm{D}^b(X), \{\mathcal{L}, \mathcal{M}\})$.
By \cite{Orlov}, $\mathrm{D}^b(\mathbb{P}_X(\mathcal{O} \oplus \mathcal{L} \oplus \mathcal{M}))$ admits the following semiorthogonal decompostion.
\[\mathrm{D}^b(Z) \cong \langle p^*\mathrm{D}^b(X), p^*\mathrm{D}^b(X)\otimes \mathcal{O}(1), p^*\mathrm{D}^b(X) \otimes \mathcal{O}(2) \rangle,\]
where $\mathcal{O}(1)$ is the relative hyperplane bundle of $p:\mathbb{P}(\mathcal{O}\oplus \mathcal{L} \oplus \mathcal{M}) \to X$.

We consider the right hand side as a SOD with 2 components $A_1=p^*\mathrm{D}^b(X)$ and $A_2=\langle p^*\mathrm{D}^b(X)\otimes \mathcal{O}(1), p^*\mathrm{D}^b(X) \otimes \mathcal{O}(2) \rangle$.

By \cite[Proposition 3.8]{Orlovdgcategories}, $A_2$ is equivalent to the gluing of $p^*\mathrm{D}^b(X)\otimes \mathcal{O}(1)$ and $p^*\mathrm{D}^b(X) \otimes \mathcal{O}(2)$ along the bimodule $R$ defined by
\[
R(p^*b \otimes \mathcal{O}(2), p^*a \otimes \mathcal{O}(1)) \cong \mathrm{Hom}_{A_2}(p^*a \otimes \mathcal{O}(1), p^*b \otimes \mathcal{O}(2))
\]
for $a, b \in Ob(\mathrm{D}^b(X))$. By the adjuction and the projection formula, we have

\begin{equation*}
\begin{split}
\mathrm{Hom}_{A_2}(p^*a \otimes \mathcal{O}(1), p^*b \otimes \mathcal{O}(2)) 
& \cong \mathrm{Hom}_{A_2}(p^*a, p^*b\otimes \mathcal{O}(1)) \\
& \cong \mathrm{Hom}_{\mathrm{D}^b(X)}(a, Rp_*(p^*b\otimes \mathcal{O}(1)))\\
& \cong \mathrm{Hom}_{\mathrm{D}^b(X)}(a, b \otimes \mathcal{E}).
\end{split}
\end{equation*}

Also, $\mathrm{D}^b(\mathbb{P}(\mathcal{O}\oplus \mathcal{L} \oplus \mathcal{M}))$ is equivalent to the gluing of $A_1$ and $A_2$ along the bimodule $S$ defined by
\[
S(y, x) \cong \mathrm{Hom}_{\mathrm{D}^b(\mathbb{P}(\mathcal{O}\oplus \mathcal{L} \oplus \mathcal{M})}(x, y)
\]
for $y\in Ob(A_2)$ and $x\in Ob(A_1)$.

Recall that if a dg category $\mathcal{C}$ admits a SOD with two components $\langle T_1, T_2 \rangle$, then any object of $\mathcal{C}$ can be written as a cone of a morphism from an object of $T_1$ to an object of $T_2$. Hence we may write any object $x$ of $A_2$ as $ \mathrm{Cone}(p^*b\otimes\mathcal{O}(1) \to p^*c\otimes\mathcal{O}(2))$ for some $b,c\in Ob(\mathrm{D}^b(X))$. Also, we write $y=p^*a$. Again by the adjunction formula, we have

\begin{equation*}
\begin{split}
&S(y,x)\\
&=S(p^*a, \mathrm{Cone}(p^*b\otimes\mathcal{O}(1) \to p^*c\otimes\mathcal{O}(2)) \\
&= \mathrm{Hom}_{\mathrm{D}^b(Z)}(p^*a, \mathrm{Cone}(p^*b\otimes\mathcal{O}(1) \to p^*c\otimes\mathcal{O}(2)) \\
& \cong \mathrm{Hom}_{\mathrm{D}^b(X)}(a, Rp_*\mathrm{Cone}(p^*b\otimes\mathcal{O}(1) \to p^*c\otimes\mathcal{O}(2))\\
& \cong \mathrm{Hom}_{\mathrm{D}^b(X)}(a, \mathrm{Cone}(b\otimes\mathcal{E} \to c\otimes \mathrm{Sym}^2\mathcal{E})).\\
\end{split}
\end{equation*}

Recall that an object of $\mathcal{D}=Rep(Q_{\Sigma_{\mathbb{P}^n}}, \mathrm{D}^b(X), \{\mathcal{L}, \mathcal{M}\})$ can be represented as the Figure \ref{fig:typical} where $a,b,c$ are objects in $\mathrm{D}^b(X)$.

\begin{figure}
\resizebox{0.9\textwidth}{!}{
\begin{tikzpicture}[scale=0.9]
\newcommand*{\edgelen}{4}; \newcommand*{\vertrad}{.12}; \newcommand*{\crad}{.05}
	\newcommand*{\gspace}{.9}; \newcommand*{\angdelta}{22.5};

	\node (a) [matrix] at (-5.3,0) {
		\coordinate (tl) at (0,\edgelen); \coordinate (tr) at (\edgelen,\edgelen);
		\coordinate (bl) at (0,0); \coordinate (br) at (\edgelen,0); 
		\coordinate (tm) at (\edgelen/3,\edgelen); \coordinate (bm) at (\edgelen/3,0);
		\coordinate (mm) at (\edgelen/3,\edgelen/3);
		\coordinate (lm) at (0, \edgelen/3); \coordinate (rm) at (\edgelen, \edgelen/3);
		\coordinate (lmin) at (\edgelen/3, \edgelen*2/3); \coordinate (rmin) at (\edgelen*2/3, \edgelen/3);
		\foreach \a/\b in {tl/bl, tr/br} {\draw[asdstyle] (\a) to (\b);}
		\foreach \a/\b in {tr/tl, br/bl} {\draw[asdstyle] (\a) to (\b);}
		\foreach \a/\b in {lm/rm} {\draw[asdstyle,righthairs] (\a) to (\b);}
		\foreach \a/\b in {bm/tm, tl/br} {\draw[asdstyle,lefthairs] (\a) to (\b);}
		\foreach \c in {mm} \foreach \ang in {8,...,12}
			{\draw[blue,thick] ($(\c)+(\ang*\angdelta:\vertrad+.02)$) to (\c);}
		\foreach \c in {lmin} \foreach \ang in {2,...,8}
			{\draw[blue,thick] ($(\c)+(\ang*\angdelta:\vertrad+.02)$) to (\c);}
		\foreach \c in {rmin} \foreach \ang in {0,1,2,12,13,14,15,16}
			{\draw[blue,thick] ($(\c)+(\ang*\angdelta:\vertrad+.02)$) to (\c);}
		\node (tltri) at (\edgelen*2/9, \edgelen*7/8) {$0$};
		\node (square) at (\edgelen*1.5/9, \edgelen*1.5/9) {$0$};
		\node (brtri) at (\edgelen*8.2/9, \edgelen*1.5/9) {$0$};
		\node (hexa) at (\edgelen*2/3, \edgelen*2/3) {$0$};
		\node (btrap) at (\edgelen*1/2, \edgelen*1.5/9) {$0$};
		\node (ttrap) at (\edgelen*1.5/9, \edgelen*1/2) {$0$};
		\node (center) at (\edgelen*3/7, \edgelen*3/7) {$a$};
		\draw[genmapstyle] (tltri) to (ttrap);
		\draw[genmapstyle] (tltri) to (hexa);
		\draw[genmapstyle] (square) to (ttrap);
		\draw[genmapstyle] (square) to (btrap);
		\draw[genmapstyle] (brtri) to (btrap);
		\draw[genmapstyle] (brtri) to (hexa);
		\draw[genmapstyle] (ttrap) to (center);
		\draw[genmapstyle] (btrap) to (center);
		\draw[genmapstyle] (hexa) to (center);
	\\};

	\node (b) [matrix] at (0,0) {
		\coordinate (tl) at (0,\edgelen); \coordinate (tr) at (\edgelen,\edgelen);
		\coordinate (bl) at (0,0); \coordinate (br) at (\edgelen,0); 
		\coordinate (tm) at (\edgelen/3,\edgelen); \coordinate (bm) at (\edgelen/3,0);
		\coordinate (mm) at (\edgelen/3,\edgelen/3);
		\coordinate (lm) at (0, \edgelen/3); \coordinate (rm) at (\edgelen, \edgelen/3);
		\coordinate (lmin) at (\edgelen/3, \edgelen*2/3); \coordinate (rmin) at (\edgelen*2/3, \edgelen/3);
		\foreach \a/\b in {tl/bl, tr/br} {\draw[asdstyle] (\a) to (\b);}
		\foreach \a/\b in {tr/tl, br/bl} {\draw[asdstyle] (\a) to (\b);}
		\foreach \a/\b in {lm/rm} {\draw[asdstyle,righthairs] (\a) to (\b);}
		\foreach \a/\b in {bm/tm, tl/br} {\draw[asdstyle,lefthairs] (\a) to (\b);}
		\foreach \c in {mm} \foreach \ang in {8,...,12}
			{\draw[blue,thick] ($(\c)+(\ang*\angdelta:\vertrad+.02)$) to (\c);}
		\foreach \c in {lmin} \foreach \ang in {2,...,8}
			{\draw[blue,thick] ($(\c)+(\ang*\angdelta:\vertrad+.02)$) to (\c);}
		\foreach \c in {rmin} \foreach \ang in {0,1,2,12,13,14,15,16}
			{\draw[blue,thick] ($(\c)+(\ang*\angdelta:\vertrad+.02)$) to (\c);}
		\node (tltri) at (\edgelen*2/9, \edgelen*7/8) {$0$};
		\node (square) at (\edgelen*1.5/9, \edgelen*1.5/9) {$0$};
		\node (brtri) at (\edgelen*8.2/9, \edgelen*1.5/9) {$0$};
		\node (hexa) at (\edgelen*2/3, \edgelen*2/3) {$b\mathcal{M}$};
		\node (btrap) at (\edgelen*1/2, \edgelen*1.5/9) {$b$};
		\node (ttrap) at (\edgelen*1.5/9, \edgelen*1/2) {$b\mathcal{L}$};
		\node (center) at (\edgelen*3/7, \edgelen*3/7) {$b\mathcal{E}$};
		\draw[genmapstyle] (tltri) to (ttrap);
		\draw[genmapstyle] (tltri) to (hexa);
		\draw[genmapstyle] (square) to (ttrap);
		\draw[genmapstyle] (square) to (btrap);
		\draw[genmapstyle] (brtri) to (btrap);
		\draw[genmapstyle] (brtri) to (hexa);
		\draw[genmapstyle] (ttrap) to (center);
		\draw[genmapstyle] (btrap) to (center);
		\draw[genmapstyle] (hexa) to (center);
	\\};
	
	\node (c) [matrix] at (5.3,0) {
		\coordinate (tl) at (0,\edgelen); \coordinate (tr) at (\edgelen,\edgelen);
		\coordinate (bl) at (0,0); \coordinate (br) at (\edgelen,0); 
		\coordinate (tm) at (\edgelen/3,\edgelen); \coordinate (bm) at (\edgelen/3,0);
		\coordinate (mm) at (\edgelen/3,\edgelen/3);
		\coordinate (lm) at (0, \edgelen/3); \coordinate (rm) at (\edgelen, \edgelen/3);
		\coordinate (lmin) at (\edgelen/3, \edgelen*2/3); \coordinate (rmin) at (\edgelen*2/3, \edgelen/3);
		\foreach \a/\b in {tl/bl, tr/br} {\draw[asdstyle] (\a) to (\b);}
		\foreach \a/\b in {tr/tl, br/bl} {\draw[asdstyle] (\a) to (\b);}
		\foreach \a/\b in {lm/rm} {\draw[asdstyle,righthairs] (\a) to (\b);}
		\foreach \a/\b in {bm/tm, tl/br} {\draw[asdstyle,lefthairs] (\a) to (\b);}
		\foreach \c in {mm} \foreach \ang in {8,...,12}
			{\draw[blue,thick] ($(\c)+(\ang*\angdelta:\vertrad+.02)$) to (\c);}
		\foreach \c in {lmin} \foreach \ang in {2,...,8}
			{\draw[blue,thick] ($(\c)+(\ang*\angdelta:\vertrad+.02)$) to (\c);}
		\foreach \c in {rmin} \foreach \ang in {0,1,2,12,13,14,15,16}
			{\draw[blue,thick] ($(\c)+(\ang*\angdelta:\vertrad+.02)$) to (\c);}
		\node (tltri) at (\edgelen*2/9, \edgelen*7/8) {$c\mathcal{L}\mathcal{M}$};
		\node (square) at (\edgelen*1.5/9, \edgelen*1.5/9) {$c\mathcal{L}$};
		\node (brtri) at (\edgelen*8.2/9, \edgelen*1.5/9) {$c\mathcal{M}$};
		\node (hexa) at (\edgelen*2/3, \edgelen*2/3) {$c\mathcal{M}\mathcal{E}$};
		\node (btrap) at (\edgelen*1/2, \edgelen*1.5/9) {$c\mathcal{E}$};
		\node (ttrap) at (\edgelen*1.5/9, \edgelen*1/2) {$c\mathcal{L}\mathcal{E}$};
		\node (center) at (\edgelen*4/7, \edgelen*3/7) {$c\mathrm{Sym}^2\mathcal{E}$};
		\draw[genmapstyle] (tltri) to (ttrap);
		\draw[genmapstyle] (tltri) to (hexa);
		\draw[genmapstyle] (square) to (ttrap);
		\draw[genmapstyle] (square) to (btrap);
		\draw[genmapstyle] (brtri) to (btrap);
		\draw[genmapstyle] (brtri) to (hexa);
		\draw[genmapstyle] (ttrap) to (center);
		\draw[genmapstyle] (btrap) to (center);
		\draw[genmapstyle] (hexa) to (center);
	\\};
	\node (atxt) at (-5.3, -2.7) {$(a')$};
	\node (btxt) at (0, -2.7) {$(b')$};
	\node (ctxt) at (5.3, -2.7) {$(c')$};

	\end{tikzpicture}
}
	\caption{The components of SOD}\label{fig:components}
\end{figure}
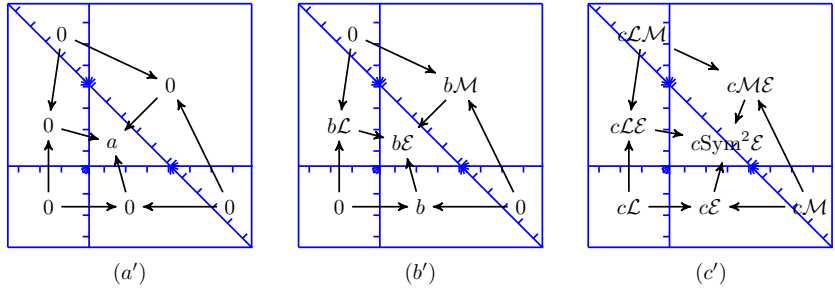

We will consider the subcategory of $Rep(Q_{\Sigma_{\mathbb{P}^n}}, \mathrm{D}^b(X), \{\mathcal{L}, \mathcal{M}\})$ generated by objects of three types described in Figure \ref{fig:components} .

Note that subcategories generated by each type of objects in Figure \ref{fig:components} are isomorphic to $\mathrm{D}^b(X)$. Let $A$, $B$, and $C$ denote the subcategories generated by objects of the form $(a')$, $(b')$, and $(c')$ in Figure  respectively.

Note that if $X$ is a point, then $A, B, C$ are exactly the subcategories of $\mathrm{Sh}_{\Lambda_{\Sigma_{\mathbb{P}^2}}}(T^2)$ generated by the objects that correspond
to $\mathcal{O}, \mathcal{O}(1), \mathcal{O}(2)$, which are exceptional collection of $\mathrm{D}^b(\mathbb{P}^2)$, under the ordinary version of CCC $\mathrm{D}^b:\mathbb{P}^2) \cong \mathrm{Sh}_{\Lambda_{\mathbb{P}^2}}(T^2)$.

We can see that they form a semiorthogonal collection. \textit{i.e.} $\mathrm{Hom}_{\mathcal{D}}(B,A), \mathrm{Hom}_{\mathcal{D}}(C,A)$, and $\mathrm{Hom}_{\mathcal{D}}(B,C)$ are all $0$.

We have $\langle B,C \rangle \cong\mathrm{D}^b(X) \times_{R'} \mathrm{D}^b(X)$ where $R'(B,C)=\mathrm{Hom}_{\mathcal{D}}(B,C)$. Any morphism from $(b')$ to $(c')$ is determined by the morphism from lower trapezoid of $(b')$ to that of $(c')$.
All other morphisms of each domain are determined by the monodromy and the universal property of the direct sum. Hence  $R'(B,C)=\mathrm{Hom}_{\mathcal{D}}(B,C)=\mathrm{Hom}_{\mathrm{D}^b(X)}(b, c\otimes\mathcal{E})$. This implies that $A_2\cong \langle B,C \rangle$ since both are the gluing of two copies of $\mathrm{D}^b(X)$ along the bimodule $\mathrm{Hom}_{\mathrm{D}^b(X)}(b, c\otimes\mathcal{E})$.
Moreover, we have 
\[
\langle A, \langle B,C\rangle \rangle=A\times_{S'} \langle B,C \rangle
\]
where $S'(x', a')=\mathrm{Hom}_{\mathcal{D}}(a',x')$ for objects $a'\in Ob(A)$ and $x\in Ob \langle B,C \rangle $. Note that $x$ can be written as $\mathrm{Cone}(b' \to c')$ for some $b'\in Ob(B)$, and $c' \in Ob(C)$.

\begin{equation*}
\begin{split}
&S'(x', a')\\
&=S'(\mathrm{Cone}(b'\to c'), a')\\
&=\mathrm{Hom}_{\mathcal{D}}(a', \mathrm{Cone}(b'\to c'))\\
&=\mathrm{Hom}_{\mathrm{D}^b(X)}(a, \mathrm{Cone}(b\otimes \mathcal{E} \to c\otimes \mathrm{Sym}^2\mathcal{E})).\\
\end{split}
\end{equation*}

Note that $A\cong \mathrm{D}^b(X) \cong A_1$. Hence both $\mathrm{D}^b(Z)=\langle A_1, A_2 \rangle $ and $\mathcal{D}=\langle A, \langle B,C\rangle \rangle$ are gluing of $\mathrm{D}^b(X)$ and  $A_2\cong \langle B,C\rangle$ along the bimodule
$$\mathrm{Hom}_{\mathrm{D}^b(X)}(a, \mathrm{Cone}(b\otimes \mathcal{E} \to c\otimes \mathrm{Sym}^2\mathcal{E})).$$
This proves that $\mathrm{D}^b(X)$ embeds into $\mathcal{D}=Rep(Q_{\Sigma_{\mathbb{P}^n}}, \mathrm{D}^b(X), \{\mathcal{L}, \mathcal{M}\})$.

To prove the theorem, now it remains to show that the subcategories $A, B, C$ generate $\mathcal{D}$.

We need to show that by taking iterative cones from $A$ to $B$ and $B$ to $C$, we generate any object in $\mathcal{D}$. Equivalently, we need to show that we obtain the zero object by taking cones from $C, B, A$ to any object in $\mathcal{D}$.

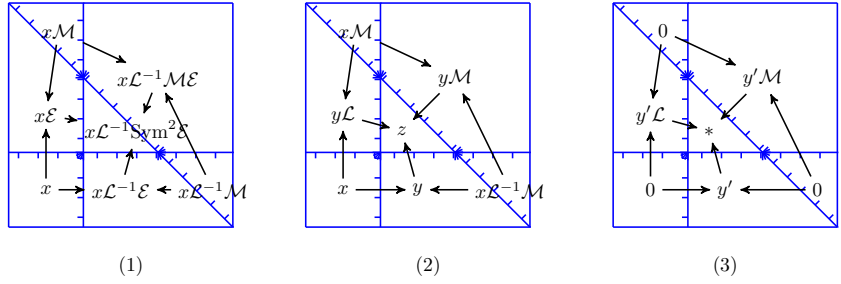
\begin{figure}
\resizebox{0.9\textwidth}{!}{
\begin{tikzpicture}
\newcommand*{\edgelen}{4}; \newcommand*{\vertrad}{.12}; \newcommand*{\crad}{.05}
	\newcommand*{\gspace}{.9}; \newcommand*{\angdelta}{22.5};

	\node (a) [matrix] at (-5.3,0) {
		\coordinate (tl) at (0,\edgelen); \coordinate (tr) at (\edgelen,\edgelen);
		\coordinate (bl) at (0,0); \coordinate (br) at (\edgelen,0); 
		\coordinate (tm) at (\edgelen/3,\edgelen); \coordinate (bm) at (\edgelen/3,0);
		\coordinate (mm) at (\edgelen/3,\edgelen/3);
		\coordinate (lm) at (0, \edgelen/3); \coordinate (rm) at (\edgelen, \edgelen/3);
		\coordinate (lmin) at (\edgelen/3, \edgelen*2/3); \coordinate (rmin) at (\edgelen*2/3, \edgelen/3);
		\foreach \a/\b in {tl/bl, tr/br} {\draw[asdstyle] (\a) to (\b);}
		\foreach \a/\b in {tr/tl, br/bl} {\draw[asdstyle] (\a) to (\b);}
		\foreach \a/\b in {lm/rm} {\draw[asdstyle,righthairs] (\a) to (\b);}
		\foreach \a/\b in {bm/tm, tl/br} {\draw[asdstyle,lefthairs] (\a) to (\b);}
		\foreach \c in {mm} \foreach \ang in {8,...,12}
			{\draw[blue,thick] ($(\c)+(\ang*\angdelta:\vertrad+.02)$) to (\c);}
		\foreach \c in {lmin} \foreach \ang in {2,...,8}
			{\draw[blue,thick] ($(\c)+(\ang*\angdelta:\vertrad+.02)$) to (\c);}
		\foreach \c in {rmin} \foreach \ang in {0,1,2,12,13,14,15,16}
			{\draw[blue,thick] ($(\c)+(\ang*\angdelta:\vertrad+.02)$) to (\c);}
		\node (tltri) at (\edgelen*2/9, \edgelen*7/8) {$x\mathcal{M}$};
		\node (square) at (\edgelen*1.5/9, \edgelen*1.5/9) {$x$};
		\node (brtri) at (\edgelen*8.2/9, \edgelen*1.5/9) {$x\mathcal{L}^{-1}\mathcal{M}$};
		\node (hexa) at (\edgelen*2/3, \edgelen*2/3) {$x\mathcal{L}^{-1}\mathcal{M}\mathcal{E}$};
		\node (btrap) at (\edgelen*1/2, \edgelen*1.5/9) {$x\mathcal{L}^{-1}\mathcal{E}$};
		\node (ttrap) at (\edgelen*1.5/9, \edgelen*1/2) {$x\mathcal{E}$};
		\node (center) at (\edgelen*4/7, \edgelen*3/7) {$x\mathcal{L}^{-1}\mathrm{Sym}^2\mathcal{E}$};
		\draw[genmapstyle] (tltri) to (ttrap);
		\draw[genmapstyle] (tltri) to (hexa);
		\draw[genmapstyle] (square) to (ttrap);
		\draw[genmapstyle] (square) to (btrap);
		\draw[genmapstyle] (brtri) to (btrap);
		\draw[genmapstyle] (brtri) to (hexa);
		\draw[genmapstyle] (ttrap) to (center);
		\draw[genmapstyle] (btrap) to (center);
		\draw[genmapstyle] (hexa) to (center);
	\\};

	\node (b) [matrix] at (0,0) {
		\coordinate (tl) at (0,\edgelen); \coordinate (tr) at (\edgelen,\edgelen);
		\coordinate (bl) at (0,0); \coordinate (br) at (\edgelen,0); 
		\coordinate (tm) at (\edgelen/3,\edgelen); \coordinate (bm) at (\edgelen/3,0);
		\coordinate (mm) at (\edgelen/3,\edgelen/3);
		\coordinate (lm) at (0, \edgelen/3); \coordinate (rm) at (\edgelen, \edgelen/3);
		\coordinate (lmin) at (\edgelen/3, \edgelen*2/3); \coordinate (rmin) at (\edgelen*2/3, \edgelen/3);
		\foreach \a/\b in {tl/bl, tr/br} {\draw[asdstyle] (\a) to (\b);}
		\foreach \a/\b in {tr/tl, br/bl} {\draw[asdstyle] (\a) to (\b);}
		\foreach \a/\b in {lm/rm} {\draw[asdstyle,righthairs] (\a) to (\b);}
		\foreach \a/\b in {bm/tm, tl/br} {\draw[asdstyle,lefthairs] (\a) to (\b);}
		\foreach \c in {mm} \foreach \ang in {8,...,12}
			{\draw[blue,thick] ($(\c)+(\ang*\angdelta:\vertrad+.02)$) to (\c);}
		\foreach \c in {lmin} \foreach \ang in {2,...,8}
			{\draw[blue,thick] ($(\c)+(\ang*\angdelta:\vertrad+.02)$) to (\c);}
		\foreach \c in {rmin} \foreach \ang in {0,1,2,12,13,14,15,16}
			{\draw[blue,thick] ($(\c)+(\ang*\angdelta:\vertrad+.02)$) to (\c);}
		\node (tltri) at (\edgelen*2/9, \edgelen*7/8) {$x\mathcal{M}$};
		\node (square) at (\edgelen*1.5/9, \edgelen*1.5/9) {$x$};
		\node (brtri) at (\edgelen*8.2/9, \edgelen*1.5/9) {$x\mathcal{L}^{-1}\mathcal{M}$};
		\node (hexa) at (\edgelen*2/3, \edgelen*2/3) {$y\mathcal{M}$};
		\node (btrap) at (\edgelen*1/2, \edgelen*1.5/9) {$y$};
		\node (ttrap) at (\edgelen*1.5/9, \edgelen*1/2) {$y\mathcal{L}$};
		\node (center) at (\edgelen*3/7, \edgelen*3/7) {$z$};
		\draw[genmapstyle] (tltri) to (ttrap);
		\draw[genmapstyle] (tltri) to (hexa);
		\draw[genmapstyle] (square) to (ttrap);
		\draw[genmapstyle] (square) to (btrap);
		\draw[genmapstyle] (brtri) to (btrap);
		\draw[genmapstyle] (brtri) to (hexa);
		\draw[genmapstyle] (ttrap) to (center);
		\draw[genmapstyle] (btrap) to (center);
		\draw[genmapstyle] (hexa) to (center);
	\\};
	
	\node (c) [matrix] at (5.3,0) {
		\coordinate (tl) at (0,\edgelen); \coordinate (tr) at (\edgelen,\edgelen);
		\coordinate (bl) at (0,0); \coordinate (br) at (\edgelen,0); 
		\coordinate (tm) at (\edgelen/3,\edgelen); \coordinate (bm) at (\edgelen/3,0);
		\coordinate (mm) at (\edgelen/3,\edgelen/3);
		\coordinate (lm) at (0, \edgelen/3); \coordinate (rm) at (\edgelen, \edgelen/3);
		\coordinate (lmin) at (\edgelen/3, \edgelen*2/3); \coordinate (rmin) at (\edgelen*2/3, \edgelen/3);
		\foreach \a/\b in {tl/bl, tr/br} {\draw[asdstyle] (\a) to (\b);}
		\foreach \a/\b in {tr/tl, br/bl} {\draw[asdstyle] (\a) to (\b);}
		\foreach \a/\b in {lm/rm} {\draw[asdstyle,righthairs] (\a) to (\b);}
		\foreach \a/\b in {bm/tm, tl/br} {\draw[asdstyle,lefthairs] (\a) to (\b);}
		\foreach \c in {mm} \foreach \ang in {8,...,12}
			{\draw[blue,thick] ($(\c)+(\ang*\angdelta:\vertrad+.02)$) to (\c);}
		\foreach \c in {lmin} \foreach \ang in {2,...,8}
			{\draw[blue,thick] ($(\c)+(\ang*\angdelta:\vertrad+.02)$) to (\c);}
		\foreach \c in {rmin} \foreach \ang in {0,1,2,12,13,14,15,16}
			{\draw[blue,thick] ($(\c)+(\ang*\angdelta:\vertrad+.02)$) to (\c);}
		\node (tltri) at (\edgelen*2/9, \edgelen*7/8) {$0$};
		\node (square) at (\edgelen*1.5/9, \edgelen*1.5/9) {$0$};
		\node (brtri) at (\edgelen*8.2/9, \edgelen*1.5/9) {$0$};
		\node (hexa) at (\edgelen*2/3, \edgelen*2/3) {$y'\mathcal{M}$};
		\node (btrap) at (\edgelen*1/2, \edgelen*1.5/9) {$y'$};
		\node (ttrap) at (\edgelen*1.5/9, \edgelen*1/2) {$y'\mathcal{L}$};
		\node (center) at (\edgelen*3/7, \edgelen*3/7) {$*$};
		\draw[genmapstyle] (tltri) to (ttrap);
		\draw[genmapstyle] (tltri) to (hexa);
		\draw[genmapstyle] (square) to (ttrap);
		\draw[genmapstyle] (square) to (btrap);
		\draw[genmapstyle] (brtri) to (btrap);
		\draw[genmapstyle] (brtri) to (hexa);
		\draw[genmapstyle] (ttrap) to (center);
		\draw[genmapstyle] (btrap) to (center);
		\draw[genmapstyle] (hexa) to (center);
	\\};
	\node (atxt) at (-5.3, -2.7) {$(1)$};
	\node (btxt) at (0, -2.7) {$(2)$};
	\node (ctxt) at (5.3, -2.7) {$(3)$};

	\end{tikzpicture}
}
	\caption{Cone from $C$}\label{fig:ConeC}
\end{figure}

Choose any object in $\mathcal{D}$, given by Figure \ref{fig:typical}((1) in Figure \ref{fig:ConeC}). Then we choose an object in $C$ given by letting $c=x\mathcal{L}^{-1}$ in $(c')$((2) in Figure \ref{fig:ConeC}). There is a natural morphism from (1) to (2) given by

\begin{align*}
& id_x: x \to x\\
& id_{x\mathcal{M}}: x\mathcal{M} \to x\mathcal{M}\\
& id_{x\mathcal{L}^{-1}\mathcal{M}}: x\mathcal{L}^{-1}\mathcal{M} \to x\mathcal{L}^{-1}\mathcal{M}
\end{align*}
and 
\[x\mathcal{L}^{-1} \mathcal{E}=x(\mathcal{L}^{-1} \oplus \mathcal{O} \oplus \mathcal{L}^{-1}\mathcal{M}) \to y  \]
is given by $0\oplus f \oplus 0$ where $f$ is a morphism from $x \to y$ in (2). Other morphisms from $x\mathcal{E} \to y\mathcal{L}$, $x\mathcal{L}^{-1}\mathcal{M}\mathcal{E} \to y\mathcal{M}$, and $x\mathcal{L}^{-1}\mathrm{Sym}^2\mathcal{E} \to z$ are defined in a similar way. 

By taking a cone of this morphism from (1) to (2), we obtain (3) in Figure \ref{fig:ConeC}, where $y'=\mathrm{Cone}(f:x\to y)$.

Then we iterate this process. We let $b=y'$ in $(b')$ in Figure \ref{fig:components} and taking cone of the canonical morphism to (3), we obtain an object in $A$. Repeating same precess again, we obtain the zero object. This completes the proof.
\end{proof}

We can extend the Theorem \ref{thm:P2} to toric $\mathbb{P}^n$-bundles for $n$ greater than $2$.
The essential ideas are the same.
In short, the idea of proof of Theorem \ref{thm:P2} is to identify the subcategories of $Rep(Q_{\Sigma_{\mathbb{P}^n}}, \mathrm{D}^b(X), \{\mathcal{L}, \mathcal{M}\})$ that are isomorphic to the  SOD components of $\mathrm{D}^b(\mathbb{P}_X(\mathcal{O} \oplus \mathcal{L} \oplus \mathcal{M}))$ and have the same gluing data.
This identification of SOD components is motivated by the ordinary CCC.

For $n$ larger than 1, 
fix a fan of $\mathbb{P}^n$ given by $\{e_1, \cdots, e_n\, -e_1-\cdots-e_n\}$ where $e_i$ denotes the $i$-th coordinate vector in $\mathbb{R}^n$. 

As in Figure \ref{fig:fltzP2}, we will perturb $\Lambda_{\Sigma_{\mathbb{P}^n}}$. We will describe an $n$-dimensional torus $T^n$ with the unit hypercube $[0,1] \times \cdots \times [0,1]$ in $\R^n$ whose parallel sides are identified. Let $\epsilon$ be a sufficiently small positive number smaller than 1. Then the perturbation $\tilde{\Lambda}_{\Sigma_{\mathbb{P}^n}}$ is obtained by pushing the $\{x_1=0\}, \cdots, \{x_n=0\}$ into the $\{x_1=\epsilon\}, \cdots, \{x_n=\epsilon\}$. In case of $n=2$, $\tilde{\Lambda}_{\Sigma_{\mathbb{P}^n}}$ is drawn in Figure \ref{fig:fltzP2}. More precisely, $\tilde{\Lambda}_{\Sigma_{\mathbb{P}^n}}$ is obtained by flowing the subset $\bigcup_{\sigma\in \Sigma }\sigma^{\perp}\times (-\sigma)$ of $\Lambda_{\Sigma_{\mathbb{P}^n}}$ under the Reeb flow for a small number of time.

The perturbation $\tilde{\Lambda}_{\Sigma_{\mathbb{P}^n}}$ divide the cube into smaller chambers.
We can classify chambers in the following way.
There is a center, which is an intersection of $(\cap_{i=1}^{n} \{\epsilon < x_i\}) \cap\{x_1+\cdots +x_n < 1\}$.
From the center, we can move to other chambers by crossing hyperplanes $\{x_1=\epsilon\}, \cdots,\{x_n=\epsilon\}$,
and $\{\sum_{i=1}^{n}x_i=1\}, \cdots \{\sum_{i=1}^{n}x_i=n-1\}$.
From the center, we can arrive at any chamber in the cube by crossing a sufficient number of hyperplanes.
Suppose one cannot cross the same hyperplane twice.
Then any chamber can be reached by crossing at most $n$ hyperplanes, and we can classify these chambers up to the number of hyperplanes to cross from the center.
This classification works in that the number of hyperplanes to cross does not depend on the choice of order of hyperplanes.
If one needs to cross $k$ hyperplanes to arrive, we say that the chamber is of $k$-step (type).
By convention, the center is the $0$-step chamber.

There is only one $0$-step chamber, the center. If $1\leq k \leq n-1$, then the number of $k$-step chambers are $\binom{n}{k}+\cdots +\binom{n}{0}$. The number of $n$-step chambers are $\binom{n}{n} +\cdots +\binom{n}{1}$. 

We can index each chamber in the following way. $(X_1, \cdots, X_n, a)$ where $X_i \in \{S, L\}$ and $a \in \{0,1,\cdots, n-1\}$. If $X_i=S$, any point in the chamber has $i$-th coordinate smaller than $\epsilon$. If $X_i=L$, then any point in the chamber has $i$-th coordinate larger than $\epsilon$. $a$ is the number of crossing slanted hyperplanes $\{\sum_{i=1}^{n}x_i=1\}, \cdots \{\sum_{i=1}^{n}x_i=n-1\}$. The center chamber is $(L, \cdots, L, 0)$. If $n=2$, for example, $(S,S,0)$ is the left bottom square, and $(L,L,1)$ is the upper right pentagon.

If $m_1, \cdots, m_i$-th entries of $(X_1, \cdots, X_n, a)$ are $S$, we will write that chamber as $(\{m_1, \cdots, m_i\}, a)$.

$\tilde{\Lambda}_{\Sigma_{\mathbb{P}^n}}$ determines a quiver $Q_{\widetilde{\Lambda}_{\Sigma_{\mathbb{P}^n}}}$. If $n=2$, the corresponding quiver is the one drawn in Figure \ref{fig:typical}.

\begin{theorem}
\label{thm:Pn}
$\mathrm{D}^b(\mathbb{P}(\mathcal{O} \oplus \mathcal{L}_1 \oplus \cdots \oplus  \mathcal{L}_n))$ is quasi-equivalent to \\
$Rep(Q_{\Sigma_{\mathbb{P}^n}}, \mathrm{D}^b(X), \{\mathcal{L}_1 \cdots , \mathcal{L}_n\})$.
\end{theorem}

\begin{proof}
By \cite{Orlov}, we know that 
\[\mathrm{D}^b(\mathbb{P}(\mathcal{O} \oplus \mathcal{L}_1 \oplus \cdots \oplus  \mathcal{L}_n))\] admits a SOD
\[
\langle p^*\mathrm{D}^b(X), p^*\mathrm{D}^b(X)\otimes \mathcal{O}(1), \cdots , p^*\mathrm{D}^b(X)\otimes \mathcal{O}(n) \rangle.
\]

We will describe subcategories of $Rep(Q_{\Sigma_{\mathbb{P}^n}}, \mathrm{D}^b(X), \{\mathcal{L}_1 \cdots , \mathcal{L}_n\})$ equivalent to each SOD component of $\mathrm{D}^b(\mathbb{P}(\mathcal{O} \oplus \mathcal{L}_1 \oplus \cdots \oplus  \mathcal{L}_n))$.

The first SOD component $p^*\mathrm{D}^b(X)$ is equivalent to the subcategory of
\[
Rep(Q_{\Sigma_{\mathbb{P}^n}}, \mathrm{D}^b(X), \{\mathcal{L}_1 \cdots , \mathcal{L}_n\})
\] consisting of following objects:

\begin{itemize}
    \item the center chamber: $a_1\in \mathrm{D}^b(X)$
    \item other chambers: $0$
\end{itemize}

The second SOD component $p^*\mathrm{D}^b(X)\otimes \mathcal{O}(1)$:
\begin{itemize}
    \item the center chamber: $a_2\mathcal{E}$
    \item $1$-step chambers: 
     \[(S, L, \cdots, L,0): a_2\]
     \[(L, S, \cdots, L,0): a_2\mathcal{L}_1\]
     \[\cdots\]
     \[(L, L, \cdots, S,0): a_2\mathcal{L}_{n-1}\]
     \[(L, L, \cdots, L, 1): a_2\mathcal{L}_n\]
    \item $2\sim n$-step chambers: $0$
\end{itemize}

The $k$-th SOD component $p^*\mathrm{D}^b(X)\otimes \mathcal{O}(k-1)$:
\begin{itemize}
    \item the center chamber: $a_k Sym^{k-1}\mathcal{E}$
    \item $1$-step chambers: 
    \[(S,L,\cdots, L,0): a_k\mathcal{L}_1Sym^{k-2}\mathcal{E}\]
    \[\cdots\]
    \[(L,\cdots,L, S,0): a_k\mathcal{L}_nSym^{k-2}\mathcal{E}\]
    \[(L, \cdots, L, 1): a_k Sym^{k-2}\mathcal{E}\]
    
    \item $i$-step chambers:
    \[\text{For }1 \leq j \leq i, (\{m_1, \cdots, m_j\},i-j): a_k\mathcal{L}_{m_1}\cdots \mathcal{L}_{m_j}Sym^{k-i-1}\mathcal{E}\]
    
    For $j=0$  
    \[(L, L, \cdots, L, i): a_kSym^{k-i-1}\mathcal{E}\]
\end{itemize}

The $n+1$-th SOD component $p^*\mathrm{D}^b(X)\otimes \mathcal{O}(n)$:

\begin{itemize}
    \item the center chamber: $Sym^n\mathcal{E}$
    \item $1$-step chambers:
    \[(S,L,\cdots,L,0): \mathcal{L}_1Sym^{n-1}\mathcal{E}\]
    \[\cdots\]
    \[(L,\cdots, L, S,0): \mathcal{L}_nSym^{n-1}\mathcal{E}\]
    \[(L, \cdots,L,1): Sym^{n-1}\mathcal{E}\]
    
    \item $k$-step chambers:
    \[(\{m_1, \cdots, m_j\}, i-j)=\mathcal{L}_{m_1}\cdots \mathcal{L}_{m_j}Sym^{k-i-1}\mathcal{E}\]
\end{itemize}

By applying exactly same steps in the proof of Theorem \ref{thm:P2}, one can prove that these subcategories, all equivalent to $\mathrm{D}^b(X)$, glue together under same bimodules and generates $Rep(Q_{\Sigma_{\mathbb{P}^n}}, \mathrm{D}^b(X), \{\mathcal{L}_1 \cdots , \mathcal{L}_n\})$. Hence $\mathrm{D}^b(\mathbb{P}(\mathcal{O} \oplus \mathcal{L}_1 \oplus \cdots \oplus  \mathcal{L}_n))$ is quasi-equivalent to $Rep(Q_{\Sigma_{\mathbb{P}^n}}, \mathrm{D}^b(X), \{\mathcal{L}_1 \cdots , \mathcal{L}_n\})$. 

\end{proof}